\documentclass[11pt]{amsart}
\usepackage
{verbatim,amssymb,amsfonts,amsthm,amsmath,amscd,ifthen,hyperref}
 \newtheorem{thm}{Theorem}[subsection]

 \newtheorem{prop}[thm]{Proposition}
 \newtheorem*{thm*}{Theorem}

\theoremstyle{definition}
\newtheorem{defn}[thm]{Definition}
\newtheorem*{defn*}{Definition}
\newtheorem{rem}[thm]{Remark}

\numberwithin{equation}{subsection}

\theoremstyle{definition}

\newcommand{\Tr}{\mathop{\mathrm{Tr}}}

\DeclareMathOperator{\SYT}{SYT} 
 
\DeclareMathOperator{\Diag}{diag} \DeclareMathOperator{\Ind}{ind}
\DeclareMathOperator{\Ch}{ch} \DeclareMathOperator{\Res}{res}
\DeclareMathOperator{\Ct}{c}
\newcommand{\TR}{\triangle}

\begin{document}
\title[Generalized characters of the symmetric group] {Generalized characters of the symmetric group}
\author
{Eugene Strahov}
\address{Department of Mathematics,  253-37,  Caltech, Pasadena, CA
91125} \email{strahov@caltech.edu}
\begin{abstract}
Normalized irreducible characters of the symmetric group $S(n)$
can be understood as  zonal spherical functions of the Gelfand
pair $(S(n)\times S(n),\Diag S(n))$. They form an orthogonal basis
in the space of the functions on the group $S(n)$ invariant with
respect to conjugations by $S(n)$. In this paper we consider a
different Gelfand pair connected with the symmetric group, that is
an  ``unbalanced'' Gelfand pair $(S(n)\times S(n-1),\Diag
S(n-1))$. Zonal spherical functions of this Gelfand pair  form an
orthogonal basis in a larger space of functions on $S(n)$, namely
in  the space of functions invariant with respect to conjugations
by $S(n-1)$. We refer to these zonal spherical functions as
normalized \textit{generalized} characters of $S(n)$. The main
discovery of the present paper is that these generalized
characters can be computed on the same level as the irreducible
characters of the symmetric group. The paper gives a
Murnaghan-Nakayama type rule, a Frobenius type formula, and an
analogue of the determinantal formula for the generalized
characters of $S(n)$.
\end{abstract}

 \maketitle
%%%%%%%%%%%%%%%%%%%%%%%%%%%%%%%%%%%%%%%%%%%%%%%%%%%%%%%%%%%%

\section{Introduction}
\subsection{Preliminaries and formulation of the problem}
One of the central goals of the representation theory of finite
groups is in computation of characters of irreducible
representations. When a group under considerations is the
symmetric group, $S(n)$, the irreducible characters can be
computed using either the Frobenius formula, or the determinantal
formula, or the Murnaghan-Nakayama rule (see, for example,
Macdonald \cite{macdonald}, Sagan \cite{sagan}, Stanley
\cite{stanley}).

Let $\Lambda$ denote the algebra of symmetric functions, which is
a graded algebra, isomorphic to the algebra of polynomials in the
power sums $p_1, p_2,\ldots $. If we define
$p_{\rho}=p_{\rho_1}p_{\rho_2}\ldots $ for each partition
$\rho=(\rho_1,\rho_2,\ldots )$\footnote{As in Macdonald
\cite{macdonald} we identify each partition with its Young
diagram.}, then the $p_{\rho}$ form a homogeneous basis in
$\Lambda$. Another natural homogeneous basis in $\Lambda$ is
formed by the Schur functions $s_{\lambda}$ indexed by Young
diagrams $\lambda$. The Frobenius formula is
$$
 p_{\rho}=\sum\limits_{\lambda\vdash
n}\chi^{\lambda}_{\rho}s_{\lambda},
$$
where $\chi^{\lambda}_{\rho}$ is the value of the irreducible
character $\chi^{\lambda}$ of the symmetric group $S(n)$ on the
conjugacy class in $S(n)$ indexed by the partition $\rho$ of $n$.
This formula is the key result in the classical theory of
characters of the symmetric group $S(n)$. It shows that the
character table is  the transition matrix between two bases
$p_{\rho}$ and $s_{\lambda}$ in the algebra of symmetric functions
$\Lambda$. The Frobenius formula follows from the fact that the
Schur functions $s_{\lambda}$ are images of $\chi^{\lambda}$ in
$\Lambda$ under a certain map. This map is called the
characteristic map, see \cite{macdonald} I, $\S$7. Thus, if we
denote this map by $\Ch$, we have
$$
s_{\lambda}=\Ch(\chi^{\lambda}).
$$

Another available result on irreducible characters of $S(n)$ is
the formula which represents an irreducible character,
$\chi^{\lambda}$, as an alternating sum of the induced characters
(i.e. the determinantal formula). Namely, denote by $\eta_{k}$ the
identity character of $S(k)$. If
$\lambda=(\lambda_1,\lambda_2,\ldots )$ is any partition of $n$,
let $\eta_{\lambda}$ denote
$\eta_{\lambda_1}\cdot\eta_{\lambda_2}\cdot\ldots$. Here the
multiplication, $f\cdot g$, between two characters, $f$ and $g$,
of, say, groups $S(k)$ and $S(m)$   is defined by the induction
$$
f\cdot g=\Ind_{S(k)\times S(m)}^{S(k+m)}(f\times g).
$$
With the above notation the irreducible character $\chi^{\lambda}$
is given by
$$
\chi^{\lambda}=\det\left(\eta_{\lambda_i-i+j}\right)_{1\leq
i,j\leq n}.
$$
Since $\Ch(\eta_{\lambda})=h_{\lambda}$, where
$h_{\lambda}=h_{\lambda_1}h_{\lambda_2}\ldots$ , and $h_{r}$ is
the $r$th complete symmetric function, the determinantal formula
for irreducible characters is equivalent to the Jacobi-Trudi
formula for the Schur symmetric functions,
$$
s_{\lambda}=\det\left(h_{\lambda_i-i+j}\right)_{1\leq i,j\leq n}.
$$

The Murnaghan-Nakayama rule is a recursive method to compute the
irreducible characters of the symmetric groups. It can be
formulated as follows. Let us say that a skew Young diagram is a
border strip if it is connected and does not contain any $2\times
2$ block of boxes. Suppose that $\pi\sigma$ is an element of the
symmetric group $S(n)$, where $\sigma$ is a cycle of length $j$,
and $\pi$ is a permutation of the remaining $n-j$ numbers of
cycle-type $\rho$, $\rho$ is a partition of $n-j$. The
Murnaghan-Nakayama rule says that the value of the irreducible
character of $S(n)$ parameterized by the Young diagram $\lambda$
on the element $\pi\sigma$ (i.e. of $\chi^{\lambda}(\pi\sigma)$)
is given by
$$
\chi^{\lambda}(\pi\sigma)=\underset{\nu\vdash
n-j}{\sum\limits_{\nu\subseteq\lambda}}\phi_{\lambda/\nu}\chi_{\rho}^{\nu}
$$
where $\phi_{\lambda/\nu}$ is a combinatorial coefficient. This
combinatorial coefficient is defined by the formula
$$
\phi_{\lambda/\nu}=\left\{%
\begin{array}{ll}
    (-1)^{\left\langle\lambda/\nu\right\rangle}, & \hbox{if $\lambda/\nu$ is a border strip;} \\
    0, & \hbox{otherwise,} \\
\end{array}%
\right.
$$
where $\left\langle\lambda/\nu\right\rangle$ is the height of a
border strip defined to be one less than the number of rows it
occupies.

The theory of characters can be reformulated in terms of Gelfand
pairs, see \cite{macdonald}, VII, $\S$1. Specifically, let $G$ be
a finite group, and $K$ be a subgroup of $G$. Denote by $C(G,K)$
the algebra of complex valued functions $f$ on $G$ (with
convolution as the multiplication) such that $f(kxk')=f(x)$ for
all $x\in G$ and $k,k'\in K$. If $C(G,K)$ is commutative, the pair
$(G,K)$ is called a \textit{Gelfand pair}, and one can associate
with $(G,K)$ the set of zonal spherical functions. Zonal spherical
functions have many remarkable properties, some of these
properties are analogous to those of group characters. In
particular, the set of zonal spherical functions defines an
orthogonal basis of $C(G,K)$, see \cite{macdonald}, VII, $\S$1.

If $K$ is a finite group, then the \textit{normalized} irreducible
characters of $K$ are closely connected with the zonal spherical
functions of the Gelfand pair $(K\times K,\Diag K)$, see
\cite{macdonald}, VII, $\S$1, Ex.9, and Section \ref{SECTION1}
below. (Here $\Diag K=\left\{(x,x): x\in K\right\}$ is the
diagonal subgroup of $K\times K$.) Explicitly, let $\chi_i$
$(1\leq i\leq r)$ be the irreducible characters of $K$, and
$\omega_i$ $(1\leq i\leq r)$ be the zonal spherical functions of
the Gelfand pair $(K\times K,\Diag K)$. Then for all elements
$x,y$ of the group $K$ the following formula holds
\begin{equation}\label{equationwcharacter}
\omega_i(x,y)=\frac{\chi_i(xy^{-1})}{\chi_i(e)}.
\end{equation}
In this situation $C(K\times K,\Diag K)$ can be identified with
the algebra of central functions $f$ defined on the group $K$,
i.e. it consists of the functions $f$ such that $f(xyx^{-1})=f(y)$
for all $x,y\in K$. In particular, if $K=S(n)$, where $S(n)$
denotes the symmetric group of symbols $1,2,\ldots ,n$, then
$C(S(n)\times S(n),\Diag S(n))$ can be identified with the algebra
of central functions defined on $S(n)$. We will denote this
algebra by $C(n)$. Since the zonal spherical functions of the
Gelfand pair $(S(n)\times S(n),\Diag S(n))$ can be expressed in
terms of  the normalized irreducible characters of $S(n)$, the
irreducible characters of $S(n)$ form an orthogonal basis in
$C(n)$, i.e. in the space of central functions.

The theory of characters in such formulation can be extended to
the Gelfand pairs $\left(S(2n),H(n)\right)$ where $H(n)$ is the
hyperoctahedral group of degree $n$ , see \cite{macdonald}, VII,
$\S$2. In particular, the zonal polynomials (i.e. the Jack
symmetric functions $J_{\lambda}^{(\alpha)}$ with parameter
$\alpha=2$) are the images of the zonal spherical functions of
$(S(2n),H(n))$ under a characteristic map, see \cite{macdonald}
VII, $\S$2.

However, there exist ``unbalanced" Gelfand pairs,
$\left(S(n)\times S(n-1),\Diag S(n-1)\right)$. (For a proof that
$\left(S(n)\times S(n-1),\Diag S(n-1)\right)$ is a Gelfand pair
see Travis \cite{travis}, Brender \cite{bender}, and Section 2
below). The algebra $C(S(n)\times S(n-1),\Diag S(n-1))$ can be
identified with the algebra of complex valued functions defined on
the group $S(n)$ invariant with respect to the conjugations by the
subgroup $S(n-1)$ of $S(n)$.  It consists of the functions $f$ on
$S(n)$ such that $f(xyx^{-1})=f(y)$ for all $x\in S(n-1)$ and
$y\in S(n)$. We will denote this algebra by $C'(n)$.

Thus, instead of conventional conjugacy classes in the symmetric
group $S(n)$ we will deal with the conjugacy classes in $S(n)$
defined with respect to $S(n-1)$. Suppose that the symmetric group
$S(n)$ is realized as the group of permutations of the set
$\{1,2,\ldots ,n\}$, and the subgroup $S(n-1)$ of $S(n)$ is
realized as the group of permutations of the set $\{1,2,\ldots
,n-1\}$. Then the conjugacy classes with respect to $S(n-1)$ can
be parameterized by pairs of Young diagrams
$(\rho,\sigma\nearrow\rho)$ where $\rho$ is a Young diagram with
$n$ boxes, $\sigma$ is a Young diagram with $n-1$ boxes, and the
notation $\sigma\nearrow\rho$ means that $\rho$ is obtained from
$\sigma$ by adding one box, see Section \ref{SectionTableaux} for
details. It turns out that the zonal spherical functions of
$(S(n)\times S(n-1),\Diag S(n-1))$ can be parameterized by pairs
of Young diagrams $(\lambda,\mu\nearrow\lambda)$ as well. Let us
denote these zonal spherical functions by
$\omega^{\lambda,\mu\nearrow\lambda}$. For all elements $x\in
S(n)$ and $y\in S(n-1)$ these spherical functions can be expressed
explicitly as:
$$
\omega^{\lambda,\mu\nearrow\lambda}(x,y)=\frac{1}{(n-1)!}\sum\limits_{z\in
S(n-1)}\chi^{\lambda}(xz)\chi^{\mu}(yz),
$$
see Brender \cite{bender} and Section 2 below for more details.

Motivated by relation (\ref{equationwcharacter}) (which connects
irreducible characters of a finite group $K$ with zonal spherical
functions of ``balanced" Gelfand pairs $(K\times K,\Diag K))$ we
define a \textit{generalized character} of $S(n)$,
$\Gamma^{\lambda,\mu\nearrow\lambda}$, by the following formula:
$$
\Gamma^{\lambda,\mu\nearrow\lambda}(x)=\chi^{\mu}(e)\cdot
\omega^{\lambda,\mu\nearrow\lambda}(x,e).
$$
Note that functions $\Gamma^{\lambda,\mu\nearrow\lambda}$ on the
group $S(n)$ introduced above define an orthogonal basis in the
space $C'(n)$, in contrast to the irreducible characters of the
group $S(n)$ which define an orthogonal basis in the subspace
$C(n)$ of $C'(n)$.  For this reason we refer to
$\Gamma^{\lambda,\mu\nearrow\lambda}$ as to \textit{generalized
characters } of $S(n)$.

   The present
paper aims to compute these generalized characters on the same
level as it is done in the classical theory for the irreducible
characters of the finite symmetric group. We find analogues of the
Frobenius formula, of the determinantal formula and of the
Murnaghan-Nakayama rule for these generalized characters. Thus we
find the same type of results for these objects as for the
irreducible characters of the symmetric group.

\subsection{Summary of results}
\subsubsection{A Murnaghan-Nakayama type rule} For the generalized
characters $\Gamma^{\lambda,\mu\nearrow\lambda}$ of the finite
symmetric group we derive a Murnaghan-Nakayma type rule, see
Theorem \ref{MAINTHEOREM}. This is the first main result of the
present paper. Its proof is based on a rational function identity,
Theorem \ref{FORMULAFORXLM},  for certain algebraic expresseions
parameterized by pairs $(\lambda,\mu\nearrow\lambda)$ of Young
diagrams.

\subsubsection{A characteristic map}
In Section \ref{SECTION5} we introduce the space $C'(n)$ of
(complex-valued) functions on $S(n)$ invariant with respect to
conjugations by elements of the subgroup $S(n-1)$ of the group
$S(n)$, and set
$$
C'=\underset{n\geq 1}{\bigoplus}\;C'(n).
$$
Then we construct an isomorphism between $C'$ and $\Lambda[t]$,
where $\Lambda[t]$ is the space of polynomials in $t$ whose
coefficients are  elements of $\Lambda$, i.e. these coefficients
are symmetric functions. This isomorphism is an analogue of the
characteristic map in the classical theory. We then define in
Section \ref{SECTION6} \textit{generalized Schur functions} as
images of the generalized characters under this map.

\subsubsection{A Frobenius type formula for the generalized characters}
The second main result of the present paper is the Frobenius type
formula for the generalized characters
$\Gamma^{\lambda,\mu\nearrow\lambda}$, equation (\ref{34}). This
formula presents the table of the generalized characters as the
transition matrix between two bases in $\Lambda[t]$. These bases
are analogues of two bases in $\Lambda$ involved in the classical
Frobenius formula, namely of the Schur symmetric functions, and of
the power sum symmetric functions.

\subsubsection{Analogues of the Jacobi-Trudi formula and of the determinantal formula}
The third main result of the present paper is in obtaining an
analogue of the Jacobi-Trudi formula for  \textit{generalized
Schur functions}, equation (\ref{655}). This formula can be also
understood as an analogue of the determinantal formula for the
irreducible characters of the symmetric group.

\subsection{Remarks  on
related works}
\subsubsection{} Zonal spherical functions on
finite groups were  studied by many authors, see e.g. Travis
\cite{travis}, Gallagher \cite{gallagher} and references therein.
The zonal spherical functions associated with the Gelfand pair
$(S(n)\times S(n-1),\Diag S(n-1))$ were considered previously by
Brender \cite{bender}. This author proves a certain general
averaging theorem  and uses it to indicate how some calculations
of values can be carried out for $(S(n)\times S(n-1),\Diag
S(n-1))$. However, any explicit formula for the zonal spherical
functions of the Gelfand pair $(S(n)\times S(n-1),\Diag S(n-1))$
has not been previously known, to the best of the author's
knowledge.

\subsubsection{} Zonal spherical functions are orthogonal, and in
some cases can be expressed in terms of hypergeometric functions.
In particular, classical discrete Hahn and Krawtchouk polynomials
can be obtained in such a way. Namely, zonal spherical functions
of the Gelfand pair $(H(n),S(n))$ (where $H(n)$ is the
hyperocthachedral group) are expressed in terms of the Krawtchouk
polynomials, see Bannai and Ito \cite{bannai}. By considering
another Gelfand pair, namely $(S(n), S(n-k)\times S(k))$, one
obtains the Hahn polynomials. There are generalizations of these
results, see  Mizukawa \cite{mizukawa}, Akazawa and Mizukawa
\cite{akazawa}, Mizukawa and Tanaka \cite{mizukawa1}.

\subsubsection{} Our derivation of
the Murnaghan-Nakayama rule for the generalized characters (see
Theorem \ref{MAINTHEOREM}) starts from a projection formula, see
Section \ref{SP}. This formula was known previously (see Travis
\cite{travis}, Brender \cite{bender}, Olshanski \cite{olshanski}).
The main idea of the proof of Theorem \ref{MAINTHEOREM} is to use
Young's orthogonal representation of the symmetric group, and sum
up the corresponding diagonal elements. A similar method was used
by Greene \cite{greene} to give a combinatorial derivation of the
classical Murnaghan-Nakayama rule for the irreducible characters
of $S(n)$. The involved summation procedure in our case is based
on a rational functional identity which is an analogue of Theorem
3.3 in Greene \cite{greene}. Note that paper by Halverson and Ram
\cite{halverson} takes a similar approach to derive the characters
of the Iwahori-Hecke algebras of type $A_{n-1}$, $B_n$ and $D_n$.
 For recent Murnaghan-Nakayama type rules for Coxeter groups,
 Hecke algebras and Lie algebras see e.g. Halverson and Ram
\cite{halverson} and references therein.
 A survey paper by Roichmann  \cite{roichman} presents different combinatorial
 versions of the classical Murnaghan-Nakayama algorithm.
\subsubsection{}
In the classical theory the Frobenius formula can be viewed as a
result of the Schur-Weyl duality between  $S(n)$ and
$GL(n,\mathbb{C})$ in tensors
$$
\underset{\mbox{$n$ times}}{\underbrace{\mathbb{C}^N\otimes
\mathbb{C}^N\otimes\ldots \otimes \mathbb{C}^N}},
$$
see, for example Goodmann and Wallach \cite{goodman}. In
particular, it follows that Schur symmetric functions are
irreducible characters of $GL(n,\mathbb{C})$. It will be
interesting whether there is some analogue of the Schur-Weyl
duality in the case considered in the present paper. The paper by
Knop \cite{knop} studies zonal spherical functions for the Gelfand
pairs $(G,K)$ where $G$ is a classical group, but $G/K$ is not a
symmetric space. One of the Gelfand pairs of this type is
$(GL(n,\mathbb{C})\times GL(n-1,\mathbb{C}),\Diag
GL(n-1,\mathbb{C}))$.  The author is not aware if there exists a
relation between the spherical functions of
$(GL(n,\mathbb{C})\times GL(n-1,\mathbb{C}),\Diag
GL(n-1,\mathbb{C}))$ and the generalized Schur functions
introduced in the present paper.
 \\

 \textbf{Acknowledgement.} I am very grateful to
 Grigori Olshanski for introducing me to the problem solved in this
 paper, and
 for numerous valuable discussions.
 In particular, the way to introduce the characteristic map as in
 Section \ref{SECTION5} was suggested by Grigori Olshanski.
 I also thank
Alexei Borodin for his interest in this work, and  for his help
with computer simulations which proved to be useful in the
derivation of the Murnaghan-Nakayama type formula.
%%%%%%%%%%%%%%%%%%%%%%%%%%%%%%%%%%%%%%%%%%%%%%%%%%%%%%%%%%%%%%%%%%%
\section{Gelfand pairs and spherical functions on the symmetric
group}\label{SECTION1}
\subsection{Irreducible characters of a finite group as spherical
functions }The material of this section is standard, see
Macdonald, VII, $\S1$.
 Let $G$ be a finite group. Let
$A=\mathbb{C}[G]$ be the complex group algebra of $G$. $A$ can be
identified with the space of all complex-valued functions on $G$.
From this point of view, the multiplication in $A$ is the
convolution, $G$ acts on $A$ by the rule $(xf)(y)=f(x^{-1}y)$, the
center of $A$ consists of the central functions on $G$. Moreover,
the center of $A$ has the irreducible characters of the group $G$
as a basis. The scalar product on $A$ is
$$
\left\langle f,g\right\rangle_G=\frac{1}{|G|}\sum\limits_{x\in
G}f(x)\overline{g(x)}.
$$

Let $K$ be a subgroup of $G$. Denote by $C(G,K)$ the subalgebra of
$\mathbb{C}[G]$ which consists of the functions $f$ on $G$ such
that $f(kxk')=f(x)$ for all $x\in G$ and $k, k'\in K$. Thus
$C(G,K)$ consists of the functions constant on each double coset
$KxK$ in $G$.
\begin{thm}\label{Theorem 2.1} For a subgroup $K$ of $G$ the following
conditions are
equivalent:\\
 1. The induced representation $1_K^G$ is
multiplicity-free. \\
2. The algebra $C(G,K)$ is commutative.
\end{thm}
\begin{defn}
The pair $(G,K)$ is called a Gelfand pair if the equivalent
conditions 1, 2 of Theorem \ref{Theorem 2.1} are satisfied.
\end{defn}
%%%%%%%%%%%%%%%%%%%%%%%%%%%%%%%%%%%%%%%%%%

Assume from now on that $(G,K)$ is a Gelfand pair. Then the
induced representation $1_G^K$ is a direct sum of non-isomorphic
irreducible $G$-modules,
\begin{equation}
1_K^G=\bigoplus\limits_{i=1}^sT_i.
\end{equation}
\begin{defn}
Let $\chi_i$ be the character of $T_i$. The functions $\omega_i$
defined for $x\in G$ by
\begin{equation}\label{GeqSphFunctions}
\omega_i(x)=\frac{1}{|K|}\sum\limits_{k\in K}\chi_i(x^{-1}k)
\end{equation}
are called zonal spherical functions of the Gelfand pair $(G,K)$.
\end{defn}
%%%%%%%%%%%%%%%%%%%%%%%%%%%%%%%%%%%%%%%%

\begin{prop} $(G\times G,\Diag G)$ is a Gelfand pair.
\end{prop}
\begin{proof} Let $\chi_i$ $(1\leq i\leq r)$ be the irreducible
characters of $G$. Then the irreducible characters of $G\times G$
are $\chi_i\times\chi_j$ $(1\leq i,j\leq r)$. Using the Frobenius
reciprocity we obtain
$$
\left\langle\chi_i\times\chi_j,1_{\Diag G}^{G\times
G}\right\rangle_{G\times G}=\frac{1}{|G|}\sum\limits_{x\in
G}\chi_i(x)\chi_j(x)=\left\langle\chi_i,\overline{\chi_j}
\right\rangle_G=
  \begin{cases}
    1, & i=j, \\
    0, & \text{otherwise}.
  \end{cases}
$$
This implies that $1_{\Diag G}^{G\times G}$ is multiplicity-free,
i.e. $(G\times G,\Diag G)$ is a Gelfand pair.
\end{proof}
From the proof of the Proposition above it is clear that the
character of the induced representation is
$\sum_{i=1}^r\chi_i\times \overline{\chi}_i$, and we obtain the
following expression for the spherical functions of $(G\times
G,\Diag G)$:
\begin{equation}\label{Eq.2.4.2}
\omega_i(x,y)=\frac{\chi_i(xy^{-1})}{\chi_i(e)}.
\end{equation}
This means that the normalized irreducible characters
$\frac{\chi_i(x)}{\chi_i(e)}$ of a finite group $G$ can be
understood as the zonal spherical functions $\omega_i(x,e)$ of the
Gelfand pair $(G\times G,\Diag G)$.

\subsection{Zonal spherical functions for $(S(n)\times
S(n-1),\Diag S(n-1))$-pairs} Let $S(n)$ be the group of
permutations of $\{1,2,\ldots ,n\}$. Let $S(n-1)$ be a subgroup of
$S(n)$. We agree that $S(n-1)$ is realized as the group of
permutations of $\{1,2,\ldots ,n-1\}$.
\begin{prop}
$(S(n)\times S(n-1),\Diag S(n-1))$ is a Gelfand pair.
\end{prop}
\begin{proof}
Let $\lambda$ be a partition of $n$, and $\mu$ be a partition of
$n-1$. Let $\chi^{\lambda}$ be the character of the irreducible
representation of $S(n)$ parameterized by $\lambda$, and
$\chi^{\mu}$  be the character of the irreducible representation
of $S(n-1)$ parameterized by $\mu$. Define the function
$\chi^{\lambda}\times\chi^{\mu}$ on the group $S(n)\times S(n-1)$
as follows:
$$
\chi^{\lambda}\times\chi^{\mu}:\;\;(x,y)\rightarrow\chi^{\lambda}(x)
\chi^{\mu}(y)
$$
where $x\in S(n)$ and $y\in S(n-1)$. It is clear that
$\chi^{\lambda}\times\chi^{\mu}$ is an irreducible character of
the group $S(n)\times S(n-1)$.

The Frobenius reciprocity implies
$$
\left\langle\chi^{\lambda}\times\chi^{\mu},1_{\Diag
S(n-1)}^{S(n)\times S(n-1)}\right\rangle_{S(n)\times
S(n-1)}=\frac{1}{(n-1)!}\sum\limits_{g\in
S(n-1)}\Res_{S(n)}^{S(n-1)}(\chi^{\lambda})(g)\chi^{\mu}(g)
$$
where $\Res_{S(n)}^{S(n-1)}(\chi^{\lambda})$ denotes the
restriction of the character $\chi^{\lambda}$ of $S(n)$ to
$S(n-1)$.
 Using the relation
$$
\Res_{S(n)}^{S(n-1)}(\chi^{\lambda})=\sum\limits_{\nu\nearrow\lambda}
\chi^{\nu}
$$
where $\nu\nearrow\lambda$ means that $\lambda$ is obtained from
$\nu$ by adding one box, and exploiting the orthogonality of
irreducible characters we find
\begin{equation}\label{21}
\left\langle\chi^{\lambda}\times\chi^{\mu},1_{\Diag
S(n-1)}^{S(n)\times S(n-1)}\right\rangle_{S(n)\times S(n-1)}=
  \begin{cases}
    1, & \text{if}\; \mu\nearrow\lambda, \\
    0, & \text{otherwise}.
  \end{cases}
\end{equation}
Hence $1_{\Diag S(n-1)}^{S(n)\times S(n-1)}$ is multiplicity-free,
and $(S(n)\times S(n-1),\Diag S(n-1))$ is a Gelfand pair.
\end{proof}
It is not hard to see that the character $1_{\Diag
S(n-1)}^{S(n)\times S(n-1)}$ is equal to
$\sum\limits_{\lambda\vdash
n,\mu\nearrow\lambda}\chi^{\lambda}\times\chi^{\mu}$. Indeed, this
follows from the comparison of formula (\ref{21}) with the
following formula
$$
\left\langle\chi^{\lambda}\times\chi^{\mu},\sum\limits_{\rho\vdash
n,\nu\nearrow\rho}\chi^{\rho}\times\chi^{\nu}\right\rangle_{S(n)\times
S(n-1)}=\begin{cases}
    1, & \text{if}\; \mu\nearrow\lambda, \\
    0, & \text{otherwise}.
  \end{cases}
$$
\begin{prop}
The zonal spherical functions for $(S(n)\times S(n-1),\Diag
S(n-1))$-pairs are given by
\begin{equation}\label{Eq. 2.5.3}
\omega^{\lambda,\mu\nearrow\lambda}(x,y)=\frac{1}{(n-1)!}\sum\limits_{z\in
S(n-1)}\chi^{\lambda}(xz)\chi^{\mu}(yz).
\end{equation}
\end{prop}
\begin{proof}
The formula in the statement of the Proposition follows
immediately from general expression (\ref{GeqSphFunctions}) for
the zonal spherical functions, and from the fact that the
character $1_{\Diag S(n-1)}^{S(n)\times S(n-1)}$ is equal to
$\sum\limits_{\lambda\vdash
n,\mu\nearrow\lambda}\chi^{\lambda}\times\chi^{\mu}$.
\end{proof}
\subsection{Generalized characters}
Comparing (\ref{Eq.2.4.2}) and (\ref{Eq. 2.5.3}) it is rather
natural to introduce the notion of the generalized characters.
\begin{defn}
Let $\lambda$ be a partition of $n$, $\mu$ be a partition of
$n-1$, and $\mu\nearrow\lambda$. The functions on the group $S(n)$
defined by
\begin{equation}\label{DEFNGENERALIZEDCHARACTERS}
\Gamma^{\lambda,\mu\nearrow\lambda}(x)=\frac{\chi^{\mu}(e)}{(n-1)!}
\sum\limits_{y\in S(n-1)}\chi^{\lambda}(xy^{-1})\chi^{\mu}(y)
\end{equation}
will be called the generalized characters.
\end{defn}
\begin{rem}
As it is evident from equation (\ref{Eq. 2.5.3}) the normalized
generalized characters of $S(n)$ are related with the zonal
spherical functions for the Gelfand pair $(S(n)\times S(n-1),\Diag
S(n-1))$ as follows:
$$
\frac{\Gamma^{\lambda,\mu\nearrow\lambda}(x)}{\chi^{\mu}(e)}
=\omega^{\lambda,\mu\nearrow\lambda}(x,e).
$$
\end{rem}
\subsection{Properties of generalized characters}\label{Section24}
The following relations can be obtained immediately from
Definition \ref{DEFNGENERALIZEDCHARACTERS}, and from the basic
properties of the zonal spherical functions, equation (1.4),
Macdonald, VII, $\S 1$:
\begin{itemize}
  \item $\Gamma^{\lambda,\mu\nearrow\lambda}(e)=\chi^{\mu}(e).$
  \item $\Gamma^{\lambda,\mu\nearrow\lambda}(yxy^{-1})=
\Gamma^{\lambda,\mu\nearrow\lambda}(x)\;\;\mbox{for}\;\mbox{all}\;\;x\in
S(n),y\in S(n-1).$
  \item $\sum\limits_{y\in\;S(n)}\Gamma^{\lambda,\mu\nearrow\lambda}(xy^
{-1})\Gamma^{\rho,\nu\nearrow\rho}(y)=
\frac{n!}{\chi^{\lambda}(e)}\;\delta^{\lambda\rho}\delta^{\mu\nu}
\Gamma^{\lambda,\mu\nearrow\lambda}(y)$. \item $\left\langle
\Gamma^{\lambda,\mu\nearrow\lambda},\Gamma^{\rho,\nu\nearrow\rho}
\right\rangle_{S(n)}
=\frac{\chi^{\mu}(e)}{\chi^{\lambda}(e)}\delta^{\lambda\rho}\delta^
{\mu\nu}$.
\end{itemize}
Here the scalar product  is defined by
$$
\left\langle f,g\right\rangle_{S(n)}
=\frac{1}{n!}\sum\limits_{x\in S(n)}f(x)\overline{g(x)},
$$
for all functions $f$, $g$ on the group $S(n)$.
%%%%%%%%%%%%%%%%%%%%%%%%%%%%%%%%%%%%%%%%%
\subsection{Tables of generalized
characters}\label{SectionTableaux} We say that two permutations
$\pi$ and $\sigma$ from $S(n)$ are related by conjugation with
respect to $S(n-1)$ if $\pi=k\sigma k^{-1}$ for some permutation
$k$ from $S(n-1)$. The set of all permutations of $S(n)$ related
by conjugation with respect to $S(n-1)$ to a given $\pi$ is called
the conjugacy class of $\pi$ with respect to $S(n-1)$. It follows
that two permutations $\pi$ and $\sigma$ from $S(n)$ are in the
same conjugacy class with respect to the conjugation by $S(n-1)$
if they both have the same  cycle type, and the cycles containing
$n$ have the same length. Thus the conjugacy classes with respect
to conjugation by $S(n-1)$ can be parameterized by pairs
$(j,\rho)$ where $j$ takes values from $1$ to $n$, and $\rho$ is a
partition of $n-j$. Suppose we adopt this parameterization, and
suppose that a permutation $\pi$ is an element of the conjugacy
class $(j,\rho)$. If this permutation is written as a product of
disjoint cycles, then the cycle containing $n$ has length $j$, and
$\rho$ represents the remaining cycles. By simple combinatorial
arguments we find that the number of permutations in the conjugacy
class parameterized by $(j,\rho)$ is equal to
$\frac{(n-1)!}{z_{\rho}}$, where
$$
z_{\rho}=\prod\limits_{i\geq 1}i^{m_i}m_i!
$$
and $m_i=m_i(\rho)$ is the number of parts of $\rho$ equal to $i$.

Alternatively, $S(n-1)$-classes in $S(n)$ can be parameterized by
a partition of $n$ with a marked row. Let $\pi$ be a permutation
from $S(n)$ which has the cycle type $\rho$. Suppose the
permutation $\pi$ is written in terms of cycles, and the number
$n$ is contained in a cycle of length $k$, $1\leq k\leq n$. Then
we mark the lowest row of the partition $\rho$ which has the
length $k$. For example, the permutation $(127)(45)(36)$ belongs
to the conjugacy class with respect to conjugations by $S(6)$
which can be parameterized either by the marked partition
$(3^*,2,2)$, or by the pair $\rho=(2,2)$ and $j=3$. It is evident
that each marked partition corresponds to a pair
$\rho,\sigma\nearrow\rho$, where $\sigma$ is obtained from the
marked partition $\rho$ by removing one box from the marked row.
Thus, the marked partition $(3^*,2,2)$ can be represented by the
following pair of Young diagrams: $\rho=(3,2,2)$ and
$\sigma=(2,2,2)$. This implies a one-to-one correspondence between
$S(n-1)$-classes in $S(n)$ and pairs $(\rho,\sigma\nearrow\rho)$
of partitions.

 As it can be readily seen from
the properties of the generalized characters listed in Section
\ref{Section24}, the generalized characters are constant on
conjugacy classes with respect to the conjugation by the subgroup
$S(n-1)$. Thus in the tables of generalized characters the rows
(generalized characters) and columns (conjugacy classes) are
indexed in the same way. For the cases of $S(3), S(4)$ the tables
of normalized generalized characters can be obtained directly from
Definition \ref{DEFNGENERALIZEDCHARACTERS}.

\begin{center}
Degree 3\\

\begin{tabular}{|c|c|c|c|c|}
  \hline
  Class & $(3^*)$ & $(2,1^*)$ & $(2^*,1)$ & $(1,1,1^*)$ \\
  \hline
  Order & 2 & 1 & 2 & 1 \\
  \hline
  $\Gamma^{(3^*)}$ & 2 & 2 & 2 & 2 \\
  \hline
  $\Gamma^{(2,1^*)}$& -1 & 2 & -1 & 2 \\
  \hline
  $\Gamma^{(2^*,1)}$ & -1 & -2 & 1 & 2 \\
  \hline
  $\Gamma^{(1,1,1^*)}$ & 2 & -2 & -2 & 2 \\ \hline
\end{tabular}
\end{center}
\begin{center}
Degree 4\\

\begin{tabular}{|c|c|c|c|c|c|c|c|}
  \hline
  Class & $(4^*)$ & $(3,1^*)$ & $(3^*,1)$ & $(2,2^*)$ & $(2,1,1^*)$ &
$(2^*,1,1)$ & $(1,1,1,1^*)$ \\
\hline
  Order & 6 & 2 & 6 & 3 & 3 & 3 & 1 \\
  \hline
  $\Gamma^{(4^*)}$ & 6 & 6 & 6 & 6 & 6 & 6 & 6 \\
  \hline
  $\Gamma^{(3,1^*)}$ & -2 & 6 & -2 & -2 & 6 & -2 & 6 \\
  \hline
  $\Gamma^{(3^*,1)}$ & -2 & -3 & 1 & -2 & 0 & 4 & 6 \\
  \hline
  $\Gamma^{(2,2^*)}$ & 0 & -3 & -3 & 6 & 0 & 0 & 6 \\
  \hline
  $\Gamma^{(2,1,1^*)}$ & 2 & -3 & 1 & -2 & 0 & -4 & 6 \\
  \hline
  $\Gamma^{(2^*,1,1)}$ & 2 & 6 & -2 & -2 & -6 & 2 & 6 \\
  \hline
  $\Gamma^{(1,1,1,1^*)}$ & -6 & 6 & 6 & 6 & -6 & -6 & 6 \\ \hline
\end{tabular}
\end{center}
In the tables of the generalized characters presented above  the
order means the number of elements in the corresponding conjugacy
class.
\section{ A Murnaghan-Nakayama type rule for the generalized
characters}
\subsection{Notation}\label{SECTIONNotations}
We say that a box $b$ in a Young diagram $\lambda$ is in position
$(i,j)$  if $b$ is in row $i$ and column $j$ of $\lambda$. A Young
tableau of shape $\lambda$ is a filling of the Young diagram
$\lambda$ with the numbers $1,2,\ldots ,n$. Young tableaux in
which the numbers are increasing left to right across the rows and
increasing down the columns of $\lambda$ are called
\textit{standard}. We denote by $\SYT(\lambda)$ the set of the
standard Young tableaux of shape $\lambda$.

The \textit{content} of a box $b$ in $\lambda$ which is in
position $(i,j)$ in $\lambda$ is given by $j-i$. We write
$\Ct(b)=j-i$. Let $T$ be a Young tableau with $n$ boxes, and let
$k$ be a number which takes a value from $1$ to $n$. We denote by
$\Ct_T(k)$ the content of that box of $T$ which is occupied by
$k$. For example, if $T$ is given by

\setlength{\unitlength}{4pt}
\begin{picture}(20,20)
\put(20,0){\framebox(4,4){$7$}} \put(20,4){\framebox(4,4){$4$}}
\put(24,4){\framebox(4,4){$6$}} \put(20,8){\framebox(4,4){$2$}}
\put(24,8){\framebox(4,4){$5$}} \put(20,12){\framebox(4,4){$1$}}
\put(24,12){\framebox(4,4){$3$}} \put(28,12){\framebox(4,4){$8$}}
\end{picture}

we have $\Ct_T(5)=0$, $\Ct_T(8)=2$, and so on.

A skew Young diagram $\lambda/\nu$ is a \textit{border strip} if
it is connected and does not contain any $2\times 2$ block of
boxes. A skew Young diagram $\lambda/\nu$ is a \textit{broken
border strip} if it does not contain any $2\times 2$ block of
boxes. Therefore, a broken border strip is a union of connected
components, each of which is a border strip.

We introduce the following terminology (cf. Halverson and Ram
\cite{halverson}). A \textit{sharp corner} in a broken border
strip is a box with a box below it and  a box to its right. On the
contrary, a \textit{dull box} of a broken border strip is a box
with no box  to its right, and  no box below it. Alternatively, a
dull box can be characterized by the following property: if we
remove a dull box from a broken border strip $\lambda/\nu$ we
obtain a new broken border strip $\mu/\nu$ such that
$\mu\nearrow\lambda$.

 The following figure  shows a broken  border strip with two connected
components where each of the sharp corners has been marked with an
\textbf{s} and each of the dull boxes has been marked with a
\textbf{d}:

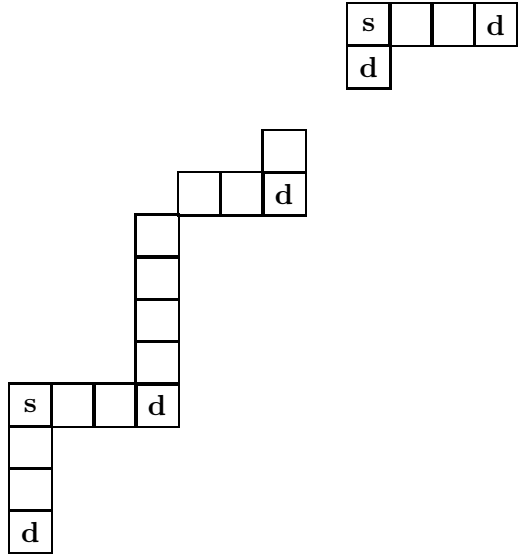
\setlength{\unitlength}{4pt}
\begin{picture}(50,60)
\put(20,0){\framebox(4,4){$\textbf{d}$}}
\put(20,4){\framebox(4,4)} \put(20,8){\framebox(4,4)}
\put(20,12){\framebox(4,4){$\textbf{s}$}}
\put(24,12){\framebox(4,4)} \put(28,12){\framebox(4,4)}
\put(32,12){\framebox(4,4){$\textbf{d}$}}
\put(32,16){\framebox(4,4)} \put(32,20){\framebox(4,4)}
\put(32,24){\framebox(4,4)} \put(32,28){\framebox(4,4)}
\put(36,32){\framebox(4,4)} \put(40,32){\framebox(4,4)}
\put(44,32){\framebox(4,4){$\textbf{d}$}}
\put(44,36){\framebox(4,4)}
\put(52,44){\framebox(4,4){$\textbf{d}$}}
\put(52,48){\framebox(4,4){$\textbf{s}$}}
\put(56,48){\framebox(4,4)} \put(60,48){\framebox(4,4)}
\put(64,48){\framebox(4,4){$\textbf{d}$}}
\end{picture}

 Recall
(see
 Macdonald \cite{macdonald}) that the
\textit{height} of a connected component (i.e. of a border strip)
is defined to be one less then the number of rows it occupies. We
define the height of a broken border strip $\lambda/\nu$ as the
sum of the heights of the connected components of $\lambda/\nu$.
We will denote the height of a broken border strip $\lambda/\nu$
by $\left\langle\lambda/\nu\right\rangle$.

\subsection{The formula for the generalized
characters}\label{S32}
\begin{thm}\label{MAINTHEOREM}(A Murnaghan-Nakayama type rule.)
Let $\lambda$ be a partition of $n$, $\mu$ be a partition of
$n-1$, and $\mu$ be obtained from $\lambda$ by removing one box.
Denote by $\Gamma^{(\lambda,\mu\nearrow\lambda)}_{(j,\rho)}$,
where $1\leq j\leq n$, and $\rho$ is a partition of $n-j$, the
value of the generalized character
$\Gamma^{(\lambda,\mu\nearrow\lambda)}$ at permutations of $S(n)$
of the cycle-type $(j,\rho)$ (with respect to conjugations by
$S(n-1)$). Then
$$
\Gamma^{(\lambda,\mu\nearrow\lambda)}_{(j,\rho)}=\underset{\nu\vdash
n-j}{\sum\limits_{\nu\subseteq\mu}}\varphi_{\mu/\nu,\lambda/\nu}\chi_{\rho}^
{\nu}
$$
where $\chi^{\nu}_{\rho}$ denotes the value of the irreducible
character of $S(n-j)$ at permutations of the group $S(n-j)$ of the
cycle type $\rho$, and $\varphi_{\mu/\nu,\lambda/\nu}$ is a
combinatorial coefficient associated with the marked skew Young
diagram $\lambda/\nu$ (the box $\lambda/\mu$ of the skew Young
diagram $\lambda/\nu$ is distinguished). This combinatorial
coefficient is defined by
$$
\varphi_{\mu/\nu,\lambda/\nu}=\left\{%
\begin{array}{ll}
    (-1)^{\left\langle\lambda/\nu\right\rangle}\prod\limits_{s\in SC}
    [\Ct(\lambda/\mu)-\Ct(s)]\underset{d\neq\lambda/\mu}{\prod\limits_
{d\in DB}}[\Ct(\lambda/\mu)-\Ct(d)]^{-1},  \\
    \;\;\;\;\;\;\;\;\;\;\;\hbox{if $\lambda/\nu$ is a broken border
strip;} \\
    \\
    0, \;\hbox{otherwise} \\
\end{array}%
\right.
$$
where $SC$ and $DB$ denote the sets of sharp corners and dull
boxes in $\lambda/\nu$, respectively, and
$\left\langle\lambda/\nu\right\rangle$ is the height of
$\lambda/\nu$.
\end{thm}
For example, let us compute the value of the generalized character
 corresponding to $\lambda=(3,2,1)$ and
$\mu=(3,2)$ at the permutation $(15)(2)(346)$. This permutation
has the following cycle type: $j=3$, $\rho=(2,1)$. The formula in
the statement of Theorem \ref{MAINTHEOREM} says that
$\Gamma^{(3,2,1),(3,2)}\left((15)(2)(346)\right)$ is a sum of two
terms since there are only two Young  diagrams $\nu$ with $6-3=3$
boxes such that $\nu\subseteq\mu$. These diagrams are $(3)$ and
$(2,1)$. Correspondingly, we need to consider contributions from
two skew Young diagrams $\lambda/\nu$,

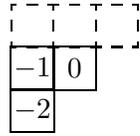
\setlength{\unitlength}{4pt}
\begin{picture}(20,20)
 \put(20,4){\framebox(4,4){$-2$}}
 \put(20,8){\framebox(4,4){$-1$}}
\put(24,8){\framebox(4,4){$0$}} \put(20,12){\dashbox(4,4)}
\put(24,12){\dashbox(4,4)} \put(28,12){\dashbox(4,4)}
\end{picture}\\
and\\
\setlength{\unitlength}{4pt}
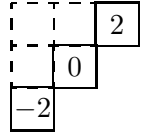
\begin{picture}(20,20)
 \put(20,4){\framebox(4,4){$-2$}}
 \put(20,8){\dashbox(4,4)}
\put(24,8){\framebox(4,4){$0$}} \put(20,12){\dashbox(4,4)}
\put(24,12){\dashbox(4,4)} \put(28,12){\framebox(4,4){$2$}}
\end{picture}\\
where the numbers in the boxes are the contents of the
corresponding boxes. We note that the height of the first diagram
is 1, and the height of the second diagram is $0$. From the
formula for the generalized characters we then find
\begin{equation}
\begin{split}
 &\Gamma^{(3,2,1),(3,2)}\left((15)(2)(346)\right)
=\\
&\chi^{(3)}_{(2,1)}(-1)^{1}\frac{-2-(-1)}{-2-0}
+\chi^{(2,1)}_{(2,1)}(-1)^{0}\frac{1}{(-2-0)(-2-2)}=-\frac{1}{2}.
\end{split}
\nonumber
\end{equation}
where we have used $\chi^{(3)}_{(2,1)}=1$,
$\chi^{(2,1))}_{(2,1)}=0$. The reader can verify that the formula
for the generalized characters stated in Theorem \ref{MAINTHEOREM}
correctly reproduces the tables of the generalized characters of
degrees 3 and 4, see Section \ref{SectionTableaux}. (Recall that
these tables were first obtained directly from the definition of
the generalized characters, equation
(\ref{DEFNGENERALIZEDCHARACTERS}).)

\section{Proof of Theorem \ref{MAINTHEOREM}}
\subsection{The projection formula}\label{SP}
A general form of a generalized character (or a spherical function
of the Gelfand pair $(S(n)\times S(n-1),\Diag S(n-1))$ is as
follows (see Travis \cite{travis}, Brender \cite{bender},
Olshanski \cite{olshanski} for more details). Suppose that
$\lambda$ and $\mu$ are any Young diagrams such that
$|\lambda|=n$, $|\mu|=n-1$, and $\mu\nearrow\lambda$ (recall that
such Young diagrams are parameters of generalized characters).
Consider the corresponding irreducible representations
$\tilde{\rho}_{\lambda}$ and $\rho_{\mu}$ of the groups $S(n)$ and
$S(n-1)$. Let $\tilde{V}^{\lambda}$ and $V^{\mu}$ be vector spaces
on which the irreducible representations $\tilde{\rho}_{\lambda}$
and $\rho_{\mu}$ are realized. Since $\rho_{\mu}$ occurs in
$\tilde{\rho}_{\lambda}|S(n)$ with multiplicity 1, we may assume
that $V^{\mu}\subset\tilde{V}^{\lambda}$. Denote by
$P^{\lambda}_{\mu}$ the projection from $\tilde{V}^{\lambda}$ onto
$V^{\mu}$. Then for any $\pi\in S(n)$
\begin{equation}\label{PROJECTIONFORMULA}
\Gamma^{\lambda,\mu\nearrow\lambda}(\pi)=\Tr\left[\tilde{\rho}_{\lambda}
(\pi)P^{\lambda}_{\mu}\right].
\end{equation}
From this projection formula it is not hard to determine the
values of the generalized characters
$\Gamma^{\lambda,\mu\nearrow\lambda}$ on $n$-cycles. This can be
done using the Jucys-Murphy elements $X_1, X_2,\ldots ,X_n$
defined by
$$
X_i=(1i)+(2i)+\ldots +(i-1i),\;\;i=1,2,\ldots ,n.
$$
In particular, $X_1=0$. (The Jucys-Murphy elements were introduced
independently in Jucys \cite{jucys} and Murphy \cite{murphy}. In
the paper by Okounkov and Vershik \cite{okounkov} the Jucys-Murphy
elements are used to give a new approach to the representation
theory of $S(n)$, see, in particular, their proof of the classical
formula for the values of irreducible characters of $S(n)$ on
$n$-cycles, Proposition 8.2 in Okounkov and Vershik
\cite{okounkov}. We also remark that in references
\cite{okounkov1,okounkov2}, \cite{olshanski,olshanski1} the
Jucys-Murphy elements are used in the infinite-dimensional
representation theory.) Now observe that $X_2X_3\ldots X_n$ is
equal to the sum of all $n$-cycles in $S(n)$. Moreover, the
eigenvalue of $X_2X_3\ldots X_n$ on any Young basis vector in
$\tilde{V}^{\lambda}$ equals
$$
(-1)^bb!(n-b-1)!
$$
if $\lambda$ is a border strip of height $b$, and equals zero
otherwise, see Okounkov and Vershik \cite{okounkov}, Proposition
8.2. Clearly, a generalized character has the same value on each
$n$-cycle, and using the projection formula, equation
(\ref{PROJECTIONFORMULA}), we find that
$\Gamma^{\lambda,\mu\nearrow\lambda}(12\ldots n))$ equals
$$
(-1)^bb!(n-b-1)!\;\frac{\dim V^{\mu}}{(n-1)!},
$$
if $\lambda$ is a border strip, and equals zero otherwise. If
$\lambda$ is a border strip of height $b$ with $n$ boxes, then
$\lambda=(a+1,b)$, and $a+b+1=n$. Note also that if $\mu$ is
obtained from $\lambda$ by removing one box, and
$\lambda=(a+1,b)$, the following two possibilities arise. The
first possibility is that the Young diagram $\mu$ has the form
$(a,b)$, and the second possibility is that the Young diagram
$\mu$ has the form $(a+1,b-1)$. In the first case $\dim
V^{\mu}=\left(\begin{array}{c}
  n-2 \\
  b
\end{array}\right)$, and in the second case $\dim V^{\mu}=
\left(\begin{array}{c}
  n-2 \\
  b-1
\end{array}\right)$.
Therefore, we obtain
\begin{equation}\label{EQGAMMAONTHECIRCLE}
\Gamma^{\lambda,\mu\nearrow\lambda}\left((1,2,\ldots ,n)\right)=
  \begin{cases}
    (-1)^{b}\frac{a}{a+b},& \lambda=(a+1,b), \mu=(a,b)\\
    (-1)^{b}\frac{b}{a+b},& \lambda=(a+1,b), \mu=(a+1,b-1)\\
    0, & \lambda\; \text{is not a border strip}.
  \end{cases}
\end{equation}
It is instructive to check that the same result follows from the
general formula in the statement of Theorem \ref{MAINTHEOREM}. If
we are interested in the value of the generalized character on the
cycle $(1,2,\ldots ,n)$ (i.e. in the value of
$\Gamma^{\lambda,\mu\nearrow\lambda}_{(j=n,\rho=\emptyset)}$),
then the sum in the expression for the generalized characters in
Theorem \ref{MAINTHEOREM} is reduced to one term only. This term
corresponds to the empty diagram, $\nu=\emptyset$, and the term
itself is equal to $\varphi_{\mu,\lambda}$. But
$\varphi_{\mu,\lambda}\neq 0$ only if $\lambda$ is a border strip.
($\lambda$ is a Young diagram, and not a skew Young diagram. Thus
$\lambda$ is necessarily connected.) Parameterizing $\lambda$ as
$(a+1,b)$ and considering two possible cases when $\mu=(a,b)$, and
$\mu=(a+1,b-1)$ we obtain the formula for the generalized
character on the cycle $(1,2,\ldots ,n)$, equation
(\ref{EQGAMMAONTHECIRCLE}).

However, the derivation of the formula for the generalized
characters in a more general situation demands an additional work.
The idea is to use Young's orthogonal representation and extend
the methods of the paper by Greene \cite{greene} to the case of
the generalized characters.
\subsection{Application of Young's orthogonal
representation}\label{Section32} Let $T$ be a Young tableau with
$n$ boxes. For each $\pi\in S(n)$, let $\pi T$ denote the Young
tableau obtained from $T$ by replacing each entry $i$ of $T$ by
its image $\pi(i)$ under $\pi$. Since $S(n)$ is generated by the
transpositions $(k,k+1)$ for $1\leq k<n$, it is sufficient for
many purposes to determine the action of these transpositions on
$T$. Now for each standard Young tableau $T$ of shape $\lambda$
with $n$ boxes, and for each transposition $(k,k+1)$ in $S(n)$ we
define
$$
\tilde{\rho}_{\lambda}\left((k,k+1)\right)T=\left\{%
\begin{array}{ll}
    A(T,k)T+B(t,k)T', & \hbox{if $T'=(k,k+1)T$ is standard;} \\
    A(T,k)T, & \hbox{otherwise} \\
\end{array}%
\right.
$$
where $A(T,k)=\left[\Ct_{T}(k+1)-\Ct_{T}(k)\right]^{-1}$ and
$B(T,k)=\sqrt{1-A^2(T,k)}$. In particular we have an action of the
group $S(n)$ on the vector space consisting of all
$\mathbb{R}$-linear combinations of the standard Young tableaux of
shape $\lambda$. Young's theorem \cite{young} says that this is a
representation of $S(n)$ which is known as Young's orthogonal
representation.

Let $\pi$ be a permutation of $S(n)$ which is in the
\textit{standard form}:
$$
\pi=(1,2,\ldots ,a_1)(a_1+1,\ldots ,a_1+a_2)\ldots
(b_{K-1}+1,\ldots ,b_K)
$$
where $\{a_1,a_2,\ldots ,a_K\}$ is the cycle-type of the
permutation $\pi$, and we have written $b_i=a_1+\ldots +a_i$ for
$i=1,2,\ldots, K$. For every Young tableau $T$ of shape $\lambda$
with $n$ boxes, and for every permutation $\pi$ of $S(n)$ which is
in the standard form we set
$$
\TR_{\pi}(T)=\underset{j\neq b_1,b_2,\ldots
,b_K}{\prod\limits_{1\leq j\leq n}}A(T,j).
$$
It turns out that the diagonal entries of the matrices in the
Young orthogonal representation can be expressed in terms of
$\TR_{\pi}(T)$. To present an explicit formula let us introduce a
linear order on the set $\SYT(\lambda)$ of the standard Young
tableaux of shape $\lambda$. Namely, if $T_i$ and $T_j$ are such
that the largest disagreeing number occurs in a lower row in
$T_j$, than $T_i$ we say that $T_i$ precedes $T_j$ in the
ordering. This is known as the \textit{last-letter ordering} of
tableaux. For example, if $\lambda=(2,2,1)$ then the last-letter
ordering is\\
\setlength{\unitlength}{4pt}
\begin{picture}(20,20)
\put(10,12){\framebox(4,4){$1$}} \put(14,12){\framebox(4,4){$4$}}
 \put(10,4){\framebox(4,4){$3$}} \put(14,8){\framebox(4,4){$5$}}
 \put(10,8){\framebox(4,4){$2$}}
 \put(22,11){{$<$}}
 %%%%%%%%%%%%%%%%%%%%%%%%%%%%%%%%%%%%%
 \put(30,12){\framebox(4,4){$1$}} \put(34,12){\framebox(4,4){$3$}}
 \put(30,4){\framebox(4,4){$4$}}  \put(34,8){\framebox(4,4){$5$}}
 \put(30,8){\framebox(4,4){$2$}}
 \put(42,11){{$<$}}
 %%%%%%%%%%%%%%%%%%%%%%%%%%%%%%%%%%%%%%%%%%%%%%%%%%%%%%%%%%%%
 \put(50,12){\framebox(4,4){$1$}} \put(54,12){\framebox(4,4){$2$}}
 \put(50,4){\framebox(4,4){$4$}}  \put(54,8){\framebox(4,4){$5$}}
 \put(50,8){\framebox(4,4){$3$}}
 \put(62,11){{$<$}}
 %%%%%%%%%%%%%%%%%%%%%%%%%%%%%%%%%%%%%%%%%%%%%%%%%%%%%%%%%%%%
 \put(70,12){\framebox(4,4){$1$}} \put(74,12){\framebox(4,4){$3$}}
 \put(70,4){\framebox(4,4){$5$}}  \put(74,8){\framebox(4,4){$4$}}
 \put(70,8){\framebox(4,4){$2$}}
 \put(82,11){{$<$}}
 %%%%%%%%%%%%%%%%%%%%%%%%%%%%%%%%%%%%%%%%%%%%%%%%%%%%%%%%%%%%
 %%%%%%%%%%%%%%%%%%%%%%%%%%%%%%%%%%%%%%%%%%%%%%%%%%%%%%%%%%%%
 \put(90,12){\framebox(4,4){$1$}} \put(94,12){\framebox(4,4){$2$}}
 \put(90,4){\framebox(4,4){$5$}}  \put(94,8){\framebox(4,4){$4$}}
 \put(90,8){\framebox(4,4){$3$}}
\end{picture}\\
With these notations we have
\begin{equation}\label{diagonalelements}
\left(\tilde\rho_{\lambda}(\pi)\right)_{i,i}=\TR_{\pi}(T_i).
\end{equation}
The proof is based on  Young's orthogonal representation
(described above in this Section). For details see  Greene
\cite{greene}, Lemma 2.4 and Ratherford \cite{rutherford}, page
43.

Formula (\ref{diagonalelements}) says that in order to compute the
irreducible character of $S(n)$ parameterized by the Young diagram
$\lambda$ we need to sum up the right-hand side of equation
(\ref{diagonalelements}) over the set $\SYT(\lambda)$ of all
standard Young tableaux of shape $\lambda$ with $n$ boxes. On the
other hand, the generalized character parameterized by the pair of
the Young diagrams $\lambda$ and $\mu$ is given by the projection
formula, equation (\ref{PROJECTIONFORMULA}). This formula implies
that the generalized characters can be represented as sums over
certain subset of $\SYT(\lambda)$. Denote this subset by
$\SYT(\lambda,\mu)$. It is clear what $\SYT(\lambda,\mu)$ is. It
consists of all standard Young tableaux of the shape $\lambda$
such that the number $n$ occupies the box $\lambda/\mu$. Therefore
we obtain
\begin{equation}\label{COMB222}
\Gamma^{(\lambda,\mu\nearrow\lambda)}(\pi)=\sum\limits_{T\in\SYT
(\lambda,\mu)}\TR_{\pi}(T).
\end{equation}
%%%%%%%%%%%%%%%%%%%%%%%%%%%%%%%%%%%%%%%%%%%%%%%%%%%%%%%%
%%%%%%%%%%%%%%%%%%%%%%%%%%%%%%%%%%%%%%%%%%%%%%%%%%%%%%%%%%%%%
%%%%%%%%%%%%%%%%%%%%%%%%%%%%%%%%%%%%%%%%%%%%%%%%%%%%%%%%%%%%
\subsection{The generalized characters as sums over sequences of
Young diagrams}\label{SEction33} Let $\lambda$ be a Young diagram
with $n$ boxes, and $\mu\nearrow\lambda$. Let $\nu$ be yet another
Young diagram such that $\nu\subseteq\mu$. Suppose that the number
of boxes in the skew Young diagram $\lambda/\nu$ equals $m$. We
set
$$
\TR(\lambda/\nu)=\sum\limits_{T\in\;\SYT(\lambda/\nu)}
\frac{1}{\left[c_T(2)-c_T(1)\right]\ldots
\left[c_T(m)-c_T(m-1)\right] }.
$$
and
$$
\TR(\mu/\nu;\lambda/\nu)=\sum\limits_{T\in\;\SYT(\mu/\nu)}
\frac{1}{\left[c_T(2)-c_T(1)\right]\ldots\left[c_T(m-1)-c_T(m-2)\right]
\left[c(\lambda/\mu)-c_T(m-1)\right] }.
$$
\begin{prop}
Let $\pi$ be a permutation of $S(n)$ which is in the standard
form, $\pi=(1,2,\ldots,a_1)(a_1+1,\ldots ,a_1+a_2)\ldots
(b_{K-1}+1,\ldots ,b_K)$, where $\{a_1,a_2,\ldots ,a_K\}$ is the
cycle type of $\pi$ and $b_i=a_1+\ldots +a_i$ for $i=1,2,\ldots
,K$. Then formula (\ref{COMB222}) for the generalized character
$\Gamma^{(\lambda,\mu\nearrow\lambda)}(\pi)$ can be rewritten as
\begin{equation}\label{COM231}
\Gamma^{\lambda,\mu\nearrow\lambda}(\pi)=
\sum\limits_{S}\TR(\lambda_1)\TR(\lambda_2/\lambda_1)\ldots\TR(\lambda_
{K-2}/\lambda_{K-1})\TR(\mu/\lambda_{K-1};\lambda/\lambda_{K-1})
\end{equation}
where the sum is over all sequences $S$ of Young diagrams
$\lambda_1$, $\lambda_2,\ldots,\lambda_{K-1}$ such that
$$
\emptyset\subseteq\lambda_1\subseteq\ldots
\subseteq\lambda_{K-1}\subseteq\mu,
$$
$\lambda_i/\lambda_{i-1}$ is a skew diagram with $a_i$ boxes for
any $i$ from $1$ to $K-1$, $\lambda_0\equiv\emptyset$, and
$\mu/\lambda_{K-1}$ is a skew diagram with $a_K-1$ boxes.
\end{prop}
\begin{proof}
As it is clear from the definition of $\TR_{\pi}(\pi)$ (see
Section \ref{Section32}) formula (\ref{COMB222}) for the
generalized character $\Gamma^{(\lambda,\mu\nearrow\lambda)}(\pi)$
can be rewritten explicitly as follows
\begin{equation}
\begin{split}
\Gamma^{(\lambda,\mu\nearrow\lambda)}(\pi)&
=\sum\limits_{T\in\SYT(\lambda,\mu)}\left[\prod\limits_{1\leq
j_1\leq b_1-1}\frac{1}{c_{T}(j_1+1)-c_T(j_1)}\right]
\left[\prod\limits_{b_1\leq j_2\leq
b_2-1}\frac{1}{c_{T}(j_2+1)-c_T(j_2)}\right]\ldots\\
&\times\left[\prod\limits_{b_{K-1}\leq j_K\leq
b_K-1}\frac{1}{c_{T}(j_K+1)-c_T(j_K)}\right].\nonumber
\end{split}
\end{equation}
In this formula we collect the terms according to positions
occupied by various segments of numbers $\{1,2,\ldots ,a_1\},
\{a_1+1,\ldots ,a_1+a_2\},\ldots ,\{b_{K-1}+1,\ldots ,n-1\}$, and
obtain
\begin{equation}
\begin{split}
\Gamma^{(\lambda,\mu\nearrow\lambda)}(\pi)&
=\sum\limits_{\emptyset\subseteq\lambda_1\subseteq\ldots\subseteq\lambda_{K-1}\subseteq\mu}
\left[\sum\limits_{T_1\in\SYT(\lambda_1)}\prod\limits_{1\leq
j_1\leq b_1-1}\frac{1}{c_{T_1}(j_1+1)-c_{T_1}(j_1)}\right]\ldots\\
&\times\left[\sum\limits_{T_{k-1}\in\SYT(\lambda_{K-1}/\lambda_{K-2})}
\prod\limits_{b_{K-2}+1\leq j_{K-1}\leq
b_{K-1}-1}\frac{1}{c_{T_{K-1}}(j_{K-1}+1)-c_{T_{K-1}}(j_{K-1})}\right]\\
&\times\left[\sum\limits_{T_{k}\in\SYT(\mu/\lambda_{K-1},\lambda/\lambda_{K-1})}
\prod\limits_{b_{K-1}+1\leq j_{K}\leq
b_{K}-1}\frac{1}{c_{T_{K}}(j_{K}+1)-c_{T_{K}}(j_K)}\right]
\nonumber
\end{split}
\end{equation}
where the sum is over Young diagrams $\lambda_1,\ldots
,\lambda_{K-1}$ with $b_1,\ldots ,b_{K-1}$ boxes correspondingly.
Now we use the fact that the products above make sense if the
tableaux (or skew tableaux) have values in any ordered set, and
these products depend only on the positions and the linear order
of symbols. This enables us to rewrite the product in the sum over
sequences of Young diagram as
$$\TR(\lambda_1)\TR(\lambda_2/\lambda_1)\ldots\TR(\lambda_
{K-2}/\lambda_{K-1})\TR(\mu/\lambda_{K-1};\lambda/\lambda_{K-1})$$
which gives the formula in the statement of the Proposition.
\end{proof}
Formula (\ref{COM231}) leads us to a Murnaghan-Nakayama type rule
for the generalized characters
$\Gamma^{\lambda,\mu\nearrow\lambda}$ provided explicit formulae
for $\TR(\lambda/\nu)$, $\TR(\mu/\nu;\lambda/\nu)$ are determined.
At this point we turn again to the work of Greene \cite{greene}.
Working in the case of irreducible characters of $S(n)$ (or
spherical functions of the ``balanced" Gelfand pair $(S(n)\times
S(n),\Diag S(n))$ Greene suggested an approach to compute
$\TR(\lambda/\nu)$. The next section reproduces these
computations.
%%%%%%%%%%%%%%%%%%%%%%%%%%%%%%%%%%%%%%%%%%%%%%%%%%%%%%%%%
%%%%%%%%%%%%%%%%%%%%%%%%%%%%%%%%%%%%%%%%%%%%%%%%%%%%%%%%%%
%%%%%%%%%%%%%%%%%%%%%%%%%%%%%%%%%%%%%%%%%%%%%%%%%%%%%%%
\subsection{Computation of $\TR(\lambda/\nu)$}\label{Section34}
In this Section $\lambda/\nu$ denotes a skew Young diagram  with
$n$ boxes. Following Greene \cite{greene}, we fix a standard
labelling of the boxes of the skew Young diagram $\lambda/\nu$. By
a standard labelling we mean the one which labels boxes in a skew
Young diagram from left to right in each row, beginning with the
first row. For example, if $\lambda=(3,2,2,1)$, $\nu=(1,1)$ the
standard labelling of the skew Young diagram $\lambda/\nu$ is\\
\setlength{\unitlength}{4pt}
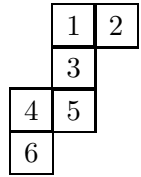
\begin{picture}(20,20)
%%%%%%%%%%%%%%%%%%%%%%%%%%%%%%%%%%%%%%%%%%%%%%%%%
\put(20,0){\framebox(4,4){$6$}} \put(20,4){\framebox(4,4){$4$}}
\put(24,4){\framebox(4,4){$5$}} \put(24,8){\framebox(4,4){$3$}}
\put(24,12){\framebox(4,4){$1$}} \put(28,12){\framebox(4,4){$2$}}
%%%%%%%%%%%%%%%%%%%%%%%%%%%%%%%%%%%%%%%%%%%%%
\end{picture}

Once the standard labelling is introduced we can regard each
standard Young tableau $T$ of shape $\lambda/\nu$ as a map
$T:\;\left\{1,2,\ldots ,n\right\}\rightarrow \left\{1,2,\ldots
,n\right\}$. This map is defined by the condition that $T(i)$ (the
value of  $T$ at an element $i$ of the set $\{1,2,\ldots ,n\}$) is
equal to the entry of the box labelled by $i$ in the standard
tableau $T$. For example, if the standard tableau $T$ is\\
\setlength{\unitlength}{4pt}
\begin{picture}(20,20)
%%%%%%%%%%%%%%%%%%%%%%%%%%%%%%%%%%%%%%%%%%%%%%%%%
\put(20,0){\framebox(4,4){$3$}} \put(20,4){\framebox(4,4){$4$}}
\put(24,4){\framebox(4,4){$6$}} \put(24,8){\framebox(4,4){$2$}}
\put(24,12){\framebox(4,4){$1$}} \put(28,12){\framebox(4,4){$5$}}
%%%%%%%%%%%%%%%%%%%%%%%%%%%%%%%%%%%%%%%%%%%%%
\end{picture}\\
then the map $T:\;\left\{1,2,3,4,5,6\right\}\rightarrow
\left\{1,2,3,4,5,6\right\}$ is defined by $T(1)=1$, $T(2)=5$,
$T(3)=2$, $T(4)=4$, $T(5)=6$, $T(6)=3$.

Let $x_1, x_2,\ldots, x_n$ be indeterminates, and assign to each
skew Young diagram with $n$ boxes a rational function constructed
with indeterminates $x_1, x_2,\ldots, x_n$. Namely, if
$\lambda/\nu$ is a skew Young diagram with $n$ boxes we assign to
$\lambda/\nu$ a rational function $X_{\lambda/\nu}(x_1,x_2,\ldots
,x_n)$ defined by
$$ X_{\lambda/\nu}(x_1,x_2,\ldots
,x_n)=\sum\limits_{T\in\;\SYT(\lambda/\nu)}
\frac{1}{\left[x_{T^{-1}(2)}-x_{T^{-1}(1)}\right]\ldots
\left[x_{T^{-1}(n)}-x_{T^{-1}(n-1)}\right] }.
$$
It is not hard to figure out how this rational function is related
with $\TR(\lambda/\nu)$. Indeed, $X_{\lambda/\nu}(x_1,x_2,\ldots
,x_n)$ is a sum of terms, and each of these terms corresponds to a
standard Young tableau of shape $\lambda/\nu$. Suppose that in
each such term we replace the indeterminates $x_1, x_2,\ldots,
x_n$ by positive integers which are contents of boxes of the Young
diagram $\lambda/\nu$. Namely, assuming the standard labelling of
$\lambda/\nu$ we replace $x_1$ by the content of the box labelled
by $1$, $x_2$ by the content of the box labelled by $2$, and so
on. Under such a replacement the rational function
$X_{\lambda/\nu}(x_1,x_2,\ldots ,x_n)$ is converted into
$\TR(\lambda/\nu)$ as it is clear from the definition of
$\TR(\lambda/\nu)$, see Section \ref{SEction33}.

The above consideration shows that the problem of computation of
$\TR(\lambda/\nu)$ can be reduced to that of finding of an
explicit formula for the rational function
$X_{\lambda/\nu}(x_1,x_2,\ldots ,x_n)$. Here we present a
rational-function identity for $X_{\lambda/\nu}(x_1,x_2,\ldots
,x_n)$. This identity is one of the central results of the paper
by Greene \cite{greene}.

Let us introduce a tableau $T_{\lambda/\nu}(x_1,\ldots ,x_n)$ of
shape $\lambda/\nu$, obtained by inserting the indeterminates
$x_1,x_2,\ldots ,x_n$ according to standard labelling of
$\lambda/\nu$. For example, if $\lambda=(3,2,2,1)$, $\nu=(1)$
then $T_{\lambda/\nu}$ is\\
\setlength{\unitlength}{4pt}
\begin{picture}(20,20)
%%%%%%%%%%%%%%%%%%%%%%%%%%%%%%%%%%%%%%%%%%%%%%%%%
\put(20,0){\framebox(4,4){$x_7$}}
\put(20,4){\framebox(4,4){$x_5$}}
\put(24,4){\framebox(4,4){$x_6$}}
\put(20,8){\framebox(4,4){$x_3$}}
\put(24,8){\framebox(4,4){$x_4$}}
\put(24,12){\framebox(4,4){$x_1$}}
\put(28,12){\framebox(4,4){$x_2$}}
%%%%%%%%%%%%%%%%%%%%%%%%%%%%%%%%%%%%%%%%%%%%%
\end{picture}\\
\begin{thm}\label{GREENETHEOREM}
If $\lambda/\nu$ is a connected skew Young diagram, the following
rational-functional identity holds
$$
X_{\lambda/\nu}(x_1,\ldots, x_n)=
\dfrac{\prod_{D_{\lambda/\nu}}(x_j-x_i)}{\prod_{R_{\lambda/\nu}}(x_j-
x_i)\prod_{C_{\lambda/\nu}}(x_j-x_i)}
$$
where $D_{\lambda}$ denotes the set of pairs $x_i$, $x_j$ with
$i<j$ which are adjacent in some diagonal in
$T_{\lambda/\nu}(x_1,\ldots ,x_n)$, $R_{\lambda/\nu}$ denotes the
set of pairs adjacent in some row of $T_{\lambda/\nu}(x_1,\ldots
,x_n)$, and $C_{\lambda/\nu}$ denotes the set of pairs adjacent in
some column of $T_{\lambda/\nu}(x_1,\ldots ,x_n)$.

If $\lambda/\nu$ is disconnected, $X_{\lambda/\nu}(x_1,\ldots
,x_n)=0$.
\end{thm}
For example,  if $\lambda=(3,2,2,1), \nu=(1)$, then Theorem
\ref{GREENETHEOREM} asserts the following expression for
$X_{\lambda/\nu}$:
\begin{equation}
X_{\lambda/\nu}=\frac{(x_6-x_3)}{(x_6-x_4)(x_4-x_1)
(x_7-x_5)(x_5-x_3)\cdot(x_2-x_1)(x_4-x_3)(x_6-x_5)}.
 \nonumber
\end{equation}
The relation between the rational function $X_{\lambda/\nu}$ and
$\TR(\lambda/\nu)$ (described above in this Section) provides us
with a formula for $\TR(\lambda/\nu)$. Namely,
\begin{equation}\label{EQUATIONTR1}
\TR(\lambda/\nu)=
  \begin{cases}
    (-1)^{<\lambda/\nu>}, & \lambda/\nu \;\text{is a border strip,}\\
    0, & \text{otherwise.} \\
  \end{cases}
\end{equation}
One of our objectives in the following study is to find an
analogue of the last written formula for
$\TR(\mu/\nu,\lambda/\nu)$ introduced previously in Section
\ref{SEction33}. If a formula for $\TR(\mu/\nu,\lambda/\nu)$ is
granted, we could  obtain an explicit expression for the
generalized character $\Gamma^{\lambda,\mu\nearrow\lambda}$, as it
is clear from equation (\ref{COM231}). In order to compute
$\TR(\mu/\nu,\lambda/\nu)$ we introduce a function
$X_{\lambda/\nu,\mu/\nu}(x_1,x_2,\ldots ,x_n)$, and derive a
rational-function identity for
$X_{\lambda/\nu,\mu/\nu}(x_1,x_2,\ldots ,x_n)$.
%%%%%%%%%%%%%%%%%%%%%%%%%%%%%%%%%%%%%%%%%%
%%%%%%%%%%%%%%%%%%%%%%%%%%%%%%%%%%%%%%%%%%%%%%%%%%%%
\subsection{ A rational-function identity related with
$\TR(\mu/\nu,\lambda/\nu)$} Let $\nu\subseteq\mu\subset\lambda$,
$\lambda/\nu$ be a skew Young diagram with $n$ boxes, $\mu/\nu$ be
a skew Young diagram with $n-1$ boxes, and $\mu\nearrow\lambda$.
As in Section \ref{Section34} we fix the standard labelling of the
boxes of $\lambda/\nu$. Suppose a box $\lambda/\mu$ is labelled by
number $\mathrm{l}(\lambda/\mu)$ under this standard labelling. We
can now regard each standard Young tableau $T$ of shape $\mu/\nu$
(filled by numbers $1,2,\ldots ,n-1$) as a map
$$
T:\; \left\{1,2,\ldots ,\mathrm{l}(\lambda/\mu)-1,
\mathrm{l}(\lambda/\mu)+1,\ldots ,n\right\}\rightarrow
\left\{1,2,\ldots ,n-1\right\},
$$
where $T(j)$ is the entry in the box labelled by $j$ in the
standard tableau $T$ of shape $\mu/\nu$.

For example, if $\lambda=(2,2,1)$, $\mu=(2,1,1)$, and $\nu=(1)$ we
fix the following labelling\\
\setlength{\unitlength}{4pt}
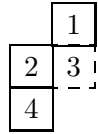
\begin{picture}(20,20)
%%%%%%%%%%%%%%%%%%%%%%%%%%%%%%%%%%%%%%%%%%%%%%%%%
\put(20,0){\framebox(4,4){$4$}} \put(20,4){\framebox(4,4){$2$}}
\put(24,4){\dashbox(4,4){$3$}} \put(24,8){\framebox(4,4){$1$}}
%%%%%%%%%%%%%%%%%%%%%%%%%%%%%%%%%%%%%%%%%%%%%
\end{picture}\\
(where the box $\lambda/\mu$ labelled by $3$ is distinguished).
Then the three standard Young tableaux of shape $\mu/\nu$, where
$\mu=(2,1,1)$ and $\nu=(1)$ can be understood as follows. The
tableau\\
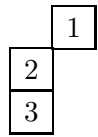
\begin{picture}(20,20)
%%%%%%%%%%%%%%%%%%%%%%%%%%%%%%%%%%%%%%%%%%%%%%%%%
\put(20,0){\framebox(4,4){$3$}} \put(20,4){\framebox(4,4){$2$}}
 \put(24,8){\framebox(4,4){$1$}}
%%%%%%%%%%%%%%%%%%%%%%%%%%%%%%%%%%%%%%%%%%%%%
\end{picture}\\
corresponds to the map $T:\; \left\{1,2,4\right\}\rightarrow
\left\{1,2,3\right\}$, under which $T(1)=1, $T(2)=2, $T(4)=3$; the
tableau
\\
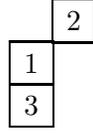
\begin{picture}(20,20)
%%%%%%%%%%%%%%%%%%%%%%%%%%%%%%%%%%%%%%%%%%%%%%%%%
\put(20,0){\framebox(4,4){$3$}} \put(20,4){\framebox(4,4){$1$}}
 \put(24,8){\framebox(4,4){$2$}}
%%%%%%%%%%%%%%%%%%%%%%%%%%%%%%%%%%%%%%%%%%%%%
\end{picture}\\
corresponds to the map $T:\; \left\{1,2,4\right\}\rightarrow
\left\{1,2,3\right\}$, defined by $T(1)=2, $T(2)=1, $T(4)=3$; and
the tableau
\\
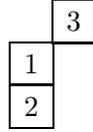
\begin{picture}(20,20)
%%%%%%%%%%%%%%%%%%%%%%%%%%%%%%%%%%%%%%%%%%%%%%%%%
\put(20,0){\framebox(4,4){$2$}} \put(20,4){\framebox(4,4){$1$}}
 \put(24,8){\framebox(4,4){$3$}}
%%%%%%%%%%%%%%%%%%%%%%%%%%%%%%%%%%%%%%%%%%%%%
\end{picture}\\
can be understood as the map $T:\; \left\{1,2,4\right\}\rightarrow
\left\{1,2,3\right\}$, under which $T(1)=3, $T(2)=1, $T(4)=2$.

With this understanding of the standard Young tableaux of shape
$\mu/\nu$ we define an algebraic expression with indeterminates
$x_1,x_2,\ldots ,x_n$. Namely, to the pair of the skew Young
diagrams $\lambda/\nu$ and $\mu/\nu$ we assign the rational
function
\begin{equation}\label{Eq351}
X_{\lambda/\nu,\mu/\nu}(x_1,\ldots
,x_n)=\sum\limits_{T\in\;\SYT(\mu/\nu)}
\frac{1}{\left[x_{T^{-1}(2)}-x_{T^{-1}(1)}\right]\ldots
\left[x_{\mathrm{l}(\lambda/\mu)}-x_{T^{-1}(n-1)}\right] }.
\end{equation}
For example, if $\lambda=(2,2,1)$, $\mu=(2,1,1)$ and $\nu=(1)$ the
right-hand side of the latest written expression takes the
following form:
$$
\frac{1}{(x_2-x_1)(x_4-x_2)(x_3-x_4)}+\frac{1}{(x_1-x_2)(x_4-x_1)(x_3-
x_4)} +\frac{1}{(x_4-x_2)(x_1-x_4)(x_3-x_1)}.
$$
We recall that  $X_{\lambda/\nu,\mu/\nu}(x_1,\ldots ,x_n)$ was
introduced with the purpose to compute $\TR(\mu/\nu,\lambda/\nu)$
included in the formula for the generalized characters
$\Gamma^{\lambda,\mu\nearrow\lambda}$, see equation
(\ref{COM231}). It follows from the definition of
$\TR(\mu/\nu,\lambda/\nu)$ that if we replace the first argument
of $X_{\lambda/\nu,\mu/\nu}$, $x_1$, by the content of the box
$\lambda/\nu$ labelled by 1, the second argument of
$X_{\lambda/\nu,\mu/\nu}$, $x_2$, by the content of the box
$\lambda/\nu$ labelled by $2$, etc, we obtain
$\TR(\mu/\nu;\lambda/\nu)$.

Now we observe a relation between $X_{\lambda/\nu}(x_1,\ldots
,x_n)$ introduced in the previous Section, and
$X_{\lambda/\nu,\mu/\nu}(x_1,\ldots ,x_n)$. Namely, we have
\begin{equation}
\begin{split}
X_{\lambda/\nu}(x_1,\ldots
,x_n)&=\sum\limits_{T\in\;\SYT(\lambda/\nu)}
\frac{1}{\left[x_{T^{-1}(2)}-x_{T^{-1}(1)}\right]\ldots
\left[x_{T^{-1}(n)}-x_{T^{-1}(n-1)}\right] }\\
&=\sum\limits_{\mu\nearrow\lambda}\sum\limits_{T\in\;\SYT(\mu/\nu)}
\frac{1}{\left[x_{T^{-1}(2)}-x_{T^{-1}(1)}\right]\ldots
\left[x_{\mathrm{l}(\lambda/\mu)}-x_{T^{-1}(n-1)}\right] }.
\nonumber
\end{split}
\end{equation}
(The last written expression is obtained if we decompose the set
$\SYT(\lambda/\nu)$ into different subsets parameterized by Young
diagrams $\mu$ such that $\mu\nearrow\lambda$. Each such subset
includes standard Young tableaux of shape $\lambda/\nu$ with $n$
boxes, and is characterized by the property that the number $n$
occupiers the same box, $\lambda/\mu$, at each tableau.) We can
also rewrite the second expression for
$X_{\lambda/\nu,\mu/\nu}(x_1,\ldots ,x_n)$ using the definition of
$X_{\lambda/\nu,\mu/\nu}(x_1,\ldots ,x_n)$, which gives
\begin{equation}\label{ZVEZDA}
X_{\lambda/\nu}(x_1,\ldots
,x_n)=\sum\limits_{\mu\nearrow\lambda}X_{\lambda/\nu,\mu/\nu}(x_1,\ldots
,x_n).
\end{equation}
Moreover, by a similar consideration we can find a recurrent
relation for $X_{\lambda/\nu,\mu/\nu}(x_1,\ldots ,x_n)$. Indeed,
\begin{equation}\label{EQRECURRENT}
\begin{split}
& X_{\lambda/\nu,\mu/\nu}(x_1,\ldots
,x_n)=\sum\limits_{T\in\;\SYT(\mu/\nu)}
\frac{1}{\left[x_{T^{-1}(2)}-x_{T^{-1}(1)}\right]\ldots
\left[x_{\mathrm{l}(\lambda/\mu)}-x_{T^{-1}(n-1)}\right] }\\
&=\sum\limits_{\gamma\nearrow\mu}\sum\limits_{T\in\;\SYT(\gamma/\nu)}
\frac{1}{\left[x_{T^{-1}(2)}-x_{T^{-1}(1)}\right]\ldots\left[x_{\mathrm{l}(\mu/
\gamma)}-x_{T^{-1}(n-2)}\right]
 \left[x_{\mathrm{l}(\lambda/\mu)}-x_{l(\mu/\gamma)}\right] }\\
&=\sum\limits_{\gamma\nearrow\mu}\frac{1}{x_{\mathrm{l}(\lambda/\mu)}-x_{\mathrm{l}(\mu/
\gamma)}}X_{\mu/\nu,\gamma/\nu}(x_1,\ldots,
\check{x}_{\mathrm{l}(\lambda/\mu)},\ldots
 ,x_n)
\end{split}
\end{equation}
where it is assumed that $\mu\supseteq\gamma\supseteq\nu$,
$\mathrm{l}(\mu/\gamma)$ is the label of the box $\mu/\gamma$
under the standard labelling, and in the last expression
$\check{x}_{\mathrm{l}(\lambda/\mu)}$ means that
$x_{\mathrm{l}(\lambda/\mu)}$ is removed from the arguments of
$X_{\mu/\nu,\gamma/\nu}$.

The above recurrent relation gives us a convenient way to
represent  $X_{\lambda/\nu}(x_1,\ldots, x_n)$ and
$X_{\mu/\nu,\lambda/\nu}(x_1,\ldots, x_n)$ graphically. For
example, we can find explicit expressions for
$X_{(3,1)}(x_1,x_2,x_3,x_4)$, $X_{(3,1),(2,1)}(x_1,x_2,x_3,x_4)$,
and $X_{(3,1),(3)}(x_1,x_2,x_3,x_4)$ using the  picture below (see
Fig. 1).
\\
 \setlength{\unitlength}{4pt}
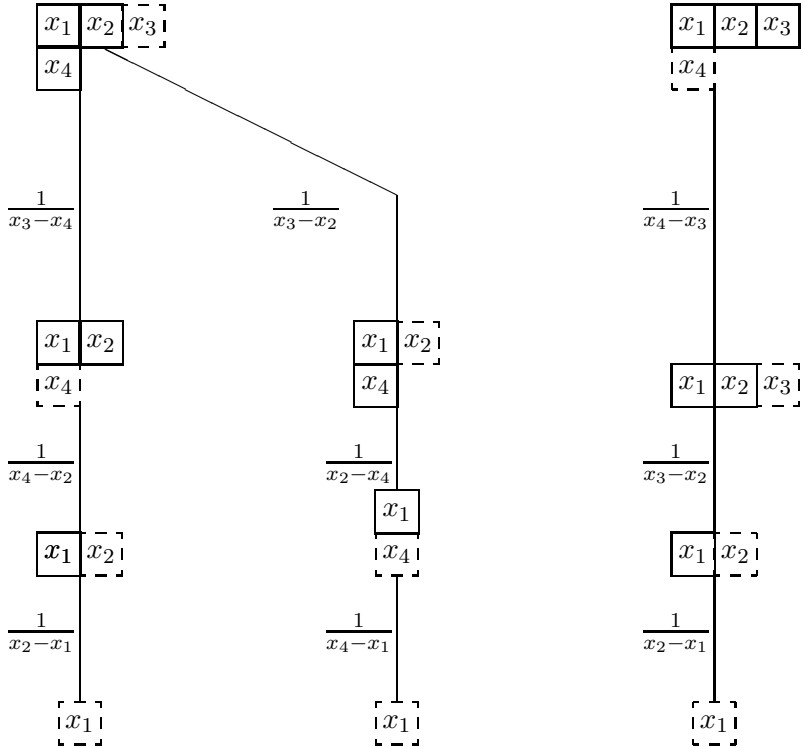
\begin{picture}(80,80)
%\put(10,-10){Fig. 1. Recurrent relations}
\put(10,0){\dashbox(4,4){$x_1$}} \put(40,0){\dashbox(4,4){$x_1$}}
\put(70,0){\dashbox(4,4){$x_1$}}
%bottom boxes
\put(5,10){{$\frac{1}{x_2-x_1}$}}
\put(5,26){{$\frac{1}{x_4-x_2}$}}
\put(5,50){{$\frac{1}{x_3-x_4}$}}
\put(30,50){{$\frac{1}{x_3-x_2}$}}
\put(65,50){{$\frac{1}{x_4-x_3}$}}
\put(65,26){{$\frac{1}{x_3-x_2}$}}
\put(65,10){{$\frac{1}{x_2-x_1}$}}
\put(35,26){{$\frac{1}{x_2-x_4}$}}
\put(35,10){{$\frac{1}{x_4-x_1}$}}
%%%%%%%%%%%
 \put(12,4){\line(0,1){12}}
\put(8,16){\framebox(4,4){$x_1$}}
\put(12,16){\dashbox(4,4){$x_2$}}
\put(8,16){\framebox(4,4){$x_1$}} \put(12,20){\line(0,1){12}}
\put(8,32){\dashbox(4,4){$x_4$}} \put(8,36){\framebox(4,4){$x_1$}}
\put(12,36){\framebox(4,4){$x_2$}}
\put(42,4){\line(0,1){12}} \put(40,16){\dashbox(4,4){$x_4$}}
\put(40,20){\framebox(4,4){$x_1$}} \put(42,24){\line(0,1){8}}
\put(38,32){\framebox(4,4){$x_4$}}
\put(38,36){\framebox(4,4){$x_1$}}
\put(42,36){\dashbox(4,4){$x_2$}}
\put(72,4){\line(0,1){12}} \put(68,16){\framebox(4,4){$x_1$}}
\put(72,16){\dashbox(4,4){$x_2$}} \put(72,20){\line(0,1){12}}
\put(76,32){\dashbox(4,4){$x_3$}}
\put(72,32){\framebox(4,4){$x_2$}}
\put(68,32){\framebox(4,4){$x_1$}}
\put(12,40){\line(0,1){22}} \put(42,40){\line(0,1){12}}
\put(72,36){\line(0,1){26}}
\put(8,62){\framebox(4,4){$x_4$}}
\put(8,66){\framebox(4,4){$x_1$}}
\put(12,66){\framebox(4,4){$x_2$}}
\put(16,66){\dashbox(4,4){$x_3$}}
\put(68,62){\dashbox(4,4){$x_4$}}
\put(68,66){\framebox(4,4){$x_1$}}
\put(72,66){\framebox(4,4){$x_2$}}
\put(76,66){\framebox(4,4){$x_3$}}
\put(42,52){\line(-2,1){28}}
\end{picture}\\
From this picture (which is just a graphical representation of the
corresponding recurrent relations, see equation
(\ref{EQRECURRENT})) we find
\begin{equation}
\begin{split}
X_{(3,1)}(x_1,x_2,x_3,x_4)=
\frac{1}{x_3-x_4}\cdot\frac{1}{x_4-x_2}\cdot\frac{1}{x_2-x_1}
&+\frac{1}{x_3-x_2}\cdot\frac{1}{x_2-x_4}\cdot\frac{1}{x_4-x_1}\\+
\frac{1}{x_4-x_3}\cdot\frac{1}{x_3-x_2}\cdot\frac{1}{x_2-x_1},
\end{split}
\nonumber
\end{equation}
\begin{equation}
\begin{split}
X_{(3,1),(2,1)}(x_1,x_2,x_3,x_4)=\frac{1}{x_3-x_4}\cdot\frac{1}{x_4-x_2}
\cdot\frac{1}{x_2-x_1}
&+\frac{1}{x_3-x_2}\cdot\frac{1}{x_2-x_4}\cdot\frac{1}{x_4-x_1},
\end{split}
\nonumber
\end{equation}
\begin{equation}
\begin{split}
X_{(3,1),(3)}(x_1,x_2,x_3,x_4)=
\frac{1}{x_4-x_3}\cdot\frac{1}{x_3-x_2}\cdot\frac{1}{x_2-x_1}.
\end{split}
\nonumber
\end{equation}
and
$X_{(3,1)}(x_1,x_2,x_3,x_4)=X_{(3,1),(2,1)}(x_1,x_2,x_3,x_4)+X_{(3,1),
(3)}(x_1,x_2,x_3,x_4)$.

Now we are in a position to give a rational-functional identity
for $X_{\lambda/\nu,\mu/\nu}(x_1,\ldots ,x_n)$. It will be clear
from the subsequent considerations that for our purposes it is
enough to consider the situation when $\lambda/\nu$ is a broken
border strip. Let $\{d_1,\ldots ,d_m\}$ be the labels of the dull
boxes of the broken border strip $\lambda/\nu$. (Recall that we
assume the standard labelling of Young diagrams introduced in the
beginning of Section \ref{Section34}). We note that one of the
numbers from the set $\{d_1,\ldots ,d_m\}$ labels the box
$\lambda/\mu$. Suppose this number is $d_i$, where the index $i$
takes some value from $1$ to $m$. Let $\{s_1,\ldots ,s_l\}$ be the
labels of sharp corners of the broken border strip $\lambda/\nu$.
Both these sets, $\{d_1,\ldots ,d_m\}$ and $\{s_1,\ldots ,s_l\}$,
are  subsets of the set $\{1,2,\ldots ,n\}$, where $n$ is the
number of boxes in $\lambda/\nu$. Recall that
$T_{\lambda/\nu}(x_1,x_2,\ldots ,x_n)$ is a tableau of shape
$\lambda/\nu$ obtained by inserting the indeterminates $x_1,\ldots
,x_n$ according to the standard labelling of $\lambda/\nu$.
\begin{thm}\label{FORMULAFORXLM}
If $\lambda/\nu$ is a broken border strip, then
$$
X_{\lambda/\nu,\mu/\nu}(x_1,\ldots ,x_n)=
\frac{\prod\limits_{j=1}^l(x_{d_i}-x_{s_j})}{\underset{j\neq i}{
\prod\limits_{j=1}^m}(x_{d_i}-x_{d_j})}\cdot\frac{1}{\prod_{R_{\lambda/
\nu}}(x_j-x_i) \prod_{C_{\lambda/\nu}}(x_j-x_i)}
$$
where $C_{\lambda/\nu}$ denotes the set of pairs $x_i,x_j$ with
$i<j$ which are adjacent in some column of
$T_{\lambda/\nu}(x_1,\ldots, x_n)$, $R_{\lambda/\nu}$ denotes the
set of pairs $x_i,x_j$ with $i<j$ which are adjacent in some row
of $T_{\lambda/\nu}(x_1,\ldots ,x_n)$, $x_{d_1},\ldots ,x_{d_m}$
are indeterminates associated with the dull boxes of
$T_{\lambda/\nu}(x_1,\ldots ,x_n)$, and $x_{s_1},\ldots ,x_{s_m}$
are indeterminates associated with the sharp boxes of
$T_{\lambda/\nu}(x_1,\ldots ,x_n)$.
\end{thm}
For example, if  the tableau $T_{\lambda/\nu}$ is that pictured
below,

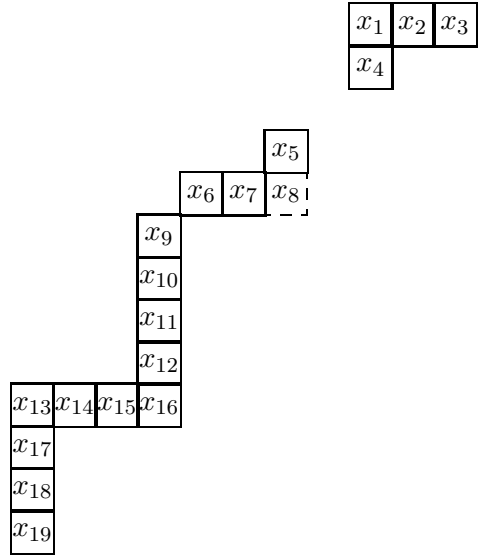
\setlength{\unitlength}{4pt}
\begin{picture}(50,60)
\put(20,0){\framebox(4,4){$x_{19}$}}
\put(20,4){\framebox(4,4){$x_{18}$}}
\put(20,8){\framebox(4,4){$x_{17}$}}
\put(20,12){\framebox(4,4){$x_{13}$}}
\put(24,12){\framebox(4,4){$x_{14}$}}
\put(28,12){\framebox(4,4){$x_{15}$}}
\put(32,12){\framebox(4,4){$x_{16}$}}
\put(32,16){\framebox(4,4){$x_{12}$}}
\put(32,20){\framebox(4,4){$x_{11}$}}
\put(32,24){\framebox(4,4){$x_{10}$}}
\put(32,28){\framebox(4,4){$x_9$}}
\put(36,32){\framebox(4,4){$x_6$}}
\put(40,32){\framebox(4,4){$x_7$}}
\put(44,32){\dashbox(4,4){$x_8$}}
\put(44,36){\framebox(4,4){$x_5$}}
\put(52,44){\framebox(4,4){$x_4$}}
\put(52,48){\framebox(4,4){$x_1$}}
\put(56,48){\framebox(4,4){$x_2$}}
\put(60,48){\framebox(4,4){$x_3$}}
\end{picture}\\
and the box corresponding to $\lambda/\mu$ is occupied by $x_8$,
Theorem \ref{FORMULAFORXLM} gives the following expression for
$X_{\lambda/\nu,\mu/\nu}$:
\begin{small}
\begin{equation}
\begin{split}
&\frac{(x_8-x_1)(x_8-x_{13})}{(x_8-x_4)(x_8-x_{16})(x_8-x_{19})}\;\cdot
\frac{1}{(x_3-x_2)(x_2-x_1)(x_8-x_7)(x_7-x_6)(x_{16}-x_{15})(x_{15}-x_{14})(x_
{14}-x_{13})}\\
&\times\frac{1}{(x_4-x_1)(x_8-x_5)(x_{16}-x_{12})(x_{12}-x_{11})
(x_{11}-x_{10})(x_{10}-x_{9})(x_{19}-x_{18})(x_{18}-x_{17})(x_{17}-x_
{13})}.
\end{split}
\nonumber
\end{equation}
\end{small}
\begin{proof}
The proof of the rational-functional identity for
$X_{\lambda/\nu,\mu/\nu}$ is by induction. In the situation when
$\lambda/\nu$ is a broken border strip with two boxes, the
asserted identity can be checked directly using the definition of
the algebraic expression $X_{\lambda/\nu,\mu/\nu}$, equation
(\ref{Eq351}). Assume that the asserted identity holds for every
broken border strip which consists of $n-1$ boxes. We are going to
show that this assumption implies the asserted identity for every
broken border strip consisting of $n$ boxes. To this end consider
all Young diagrams $\mu_1,\mu_2,\ldots ,\mu_k$ such that
$\mu_1\nearrow\mu,\mu_2\nearrow\mu,\ldots ,\mu_k\nearrow\mu$. Let
$d_1',\ldots ,d_k'$ be the labels of the boxes
$\mu/\mu_1,\mu/\mu_2,\ldots ,\mu/\mu_k$, respectively (once again,
the standard labelling of $\mu/\nu$ is assumed). With these
notations the recurrent formula for $X_{\lambda/\nu,\mu/\nu}$
derived previously in this Section can be rewritten as
\begin{equation}\label{RRR}
X_{\lambda/\nu,\mu/\nu}(x_1,\ldots,
x_n)=\sum\limits_{q=1}^k\frac{1}{x_{d_i}-x_{d'_q}}X_{\mu/\nu,\mu_q/\nu}
(x_1,\ldots, \check{x}_{d_i},\ldots ,x_n).
\end{equation}
Clearly, the boxes labelled by $d_1',d_2',\ldots ,d_k'$ are in
fact the dull boxes for the skew Young diagram $\mu/\nu$. Denote
by $s_1',s_2',\ldots ,s_p'$ the labels of the sharp corners of the
skew Young diagram $\mu/\nu$. The assumption for the case of $n-1$
boxes provides us with a rational-functional identity for
$X_{\mu/\nu,\mu_q/\nu}(x_1,\ldots, \check{x}_{d_i},\ldots ,x_n)$:
$$
X_{\mu/\nu,\mu_q/\nu}(x_1,\ldots,\check{x}_{d_i},\ldots, x_n)=
\frac{\prod\limits_{j=1}^p(x_{d'_q}-x_{s'_j})}{\underset{j\neq q}{
\prod\limits_{j=1}^k}(x_{d'_q}-x_{d'_j})}\cdot\frac{1}{\prod_{R_{\mu/
\nu}}(x_j-x_i) \prod_{C_{\mu/\nu}}(x_j-x_i)}.
$$
Insert this expression to the recurrent formula (equation
\ref{RRR})  and obtain:
$$
X_{\lambda/\nu,\mu/\nu}(x_1,\ldots,
x_n)=\left\{\sum\limits_{q=1}^k\frac{1}{x_{d_i}-x_{d'_q}}
\frac{\prod\limits_{j=1}^p(x_{d'_q}-x_{s'_j})}{\underset{j\neq q}{
\prod\limits_{j=1}^k}(x_{d'_q}-x_{d'_j})}\right\}\cdot\frac{1}{\prod_{R_
{\mu/\nu}}(x_j-x_i) \prod_{C_{\mu/\nu}}(x_j-x_i)}.
$$
In order to simplify the expression in the brackets apply the well
known algebraic identity
$$
\frac{(z-x_1)\ldots (z-x_m)}{(z-y_1)\ldots
(z-y_k)}=\sum\limits_{i=1}^k\frac{1}{z-y_i}\cdot\frac{\prod_{l=1}^m(y_i-
x_l)}{\prod^k_{j=1,j\neq i}(y_i-y_j) },
$$
which gives a rational-functional representation for
$X_{\lambda/\nu,\mu/\nu}$:
$$
X_{\lambda/\nu,\mu/\nu}(x_1,\ldots, x_n)=
\frac{\prod\limits_{j=1}^p(x_{d_i}-x_{s_j'})}{
\prod\limits_{j=1}^k(x_{d_i}-x_{d'_j})}\cdot\frac{1}{\prod_{R_{\mu/\nu}}
(x_j-x_i) \prod_{C_{\mu/\nu}}(x_j-x_i)}.
$$
It remains to show that the last written expression  actually
coincides with that in the statement of the theorem. To check
this consider two different possibilities.\\
a) The first possibility is that the broken border strips
$\lambda/\nu$ and $\mu/\nu$ have the same set of sharp corners.
Then the numerator in the last written expression is exactly the
same as that in the formula for $X_{\lambda/\nu,\mu/\nu}$ in the
statement of the theorem. If this possibility is realized it can
happen, however, that  the skew Young diagrams $\lambda/\nu$ and
$\mu/\nu$ have different dull boxes. If it is so suppose that
 the dull box of $\mu/\nu$ labelled by $d_j'$ is not a dull box of $\lambda/\nu$. Then the
corresponding term, $(x_{d_i}-x_{d'_j})$,  can be joined to one of
the products, $\prod_{R_{\mu/\nu}}(x_j-x_i)$, or
$\prod_{C_{\mu/\nu}}(x_j-x_i)$. (Note that the boxes labelled by
$d_i$ and $d'_j$ are adjacent to each other. Indeed, the box
labelled by $d_j'$ can be situated either just above the box
labelled by $d_i$ or just left to the box labelled by $d_i$. Thus
$d_i>d_j'$.) Therefore one can rewrite the denominator as
$\underset{j\neq i}
{\prod\limits_{j=1}^m}(x_{d_i}-x_{d_j})\prod_{R_{\lambda/\nu}}(x_j-x_i)
\prod_{C_{\lambda/\nu}}(x_j-x_i)$, and if  the broken border
strips $\lambda/\nu$ and $\mu/\nu$ have the same set of sharp
corners we are done.\\
b) If the second possibility is realized, i.e. the number of sharp
corners of $\mu/\nu$ is less  then that  of $\lambda/\nu$, we note
that there exists a label $s_j$ which does not coincide with any
labels $s_1',s_2',\ldots ,s_p'$. In this situation multiply the
numerator and the denominator by $(x_{d_i}-x_{s_j})$, and after
that, in the denominator, join this multiplier,
$(x_{d_i}-x_{s_j})$, to one of the  products
$\prod_{R_{\lambda/\nu}}(x_j-x_i)$ or
$\prod_{C_{\lambda/\nu}}(x_j-x_i)$. (Note that the box labelled by
$s_j$ can be situated either just above the box labelled by $d_i$,
or just  to its left. In any case $d_i>s_j$). Therefore the
situation which arises when the second possibility is realized can
be reduced to the case considered in a) (when the broken border
strips $\lambda/\nu$ and $\mu/\nu$ have the same set of sharp
corners.)
\end{proof}
\subsection{Computation of $\TR(\mu/\nu,\lambda/\nu)$}
The results of the previous Section enable us to compute
$\TR(\mu/\nu,\lambda/\nu)$, and thus to derive a formula for the
generalized characters. Let us first prove the following
\begin{prop}
Let $\lambda/\nu$ be any skew Young diagram containing a $2\times
2$ block of boxes (which means that $\lambda/\nu$ is not a broken
border strip). Then $\TR(\mu/\nu,\lambda/\nu)=0$
\end{prop}
\begin{proof}
The proof is by induction, and it is based on an application of
the following relations:
\begin{equation}\label{i}
 0=\sum\limits_{\mu\nearrow\lambda}\TR(\mu/\nu,\lambda/\nu),
\end{equation}
\begin{equation}\label{ii}
\TR(\mu/\nu,\lambda/\nu)=\sum\limits_{\gamma\nearrow\mu}\frac{1}{c
(\lambda/\mu)-c(\mu/\gamma)}. \;\TR(\gamma/\nu,\mu/\nu)
\end{equation}
(The first relation follows from the fact that
$\TR(\lambda/\nu)=0$, if $\lambda/\nu$ is not a broken border
strip, see equation (\ref{EQUATIONTR1}), and from the relation
between algebraic expressions $X_{\lambda/\nu}$ and
$X_{\lambda/\nu,\mu/\nu}$, equation (\ref{ZVEZDA}). The second
relation follows from equation (\ref{EQRECURRENT}).)

Assume that the assertion is valid for all skew Young diagram
$\lambda/\nu$ with $n-1$ boxes. Let us prove that this implies the
validity of the assertion for the case of skew Young diagrams with
$n$
boxes.\\
a) If $\lambda/\nu$  includes a $2\times 2$ block, and $\mu/\nu$
contains a $2\times 2$ block (where $\mu\nearrow\lambda$), then
$\TR(\gamma/\nu,\mu/\nu)=0$ for all $\gamma\nearrow\mu$
(assumption for $n-1$), and $\TR(\mu/\nu,\lambda/\nu)=0$ as it
follows from equation (\ref{ii}) (all terms in the sum in the
right-hand side of equation (\ref{ii})
 are zeros).\\
b) It remains to consider the case when the skew Young diagram
under consideration, $\lambda/\nu$, includes one $2\times 2$
block, but the skew Young diagram $\mu/\nu$ does not include any
$2\times 2$ blocks. Then equation (\ref{i}) gives
$$
0=\TR(\mu/\nu,\lambda/\nu)+\sum\limits_{\gamma\nearrow\lambda,
\gamma\neq\mu}\TR(\lambda/\nu,\gamma/\nu).
$$
Now, by a), all terms in the sum over Young diagrams $\gamma$ in
the right-hand side of the equation above are zeros since all
diagrams $\gamma/\nu$ must include a $2\times 2$ block of boxes
(the same $2\times 2$ block which is included in $\lambda/\nu$).
Therefore $\TR(\lambda/\nu,\mu/\nu)=0$.
\end{proof}
\begin{prop}
If $\lambda/\nu$ is any skew Young diagram,  then
$$
\TR(\mu/\nu,\lambda/\nu)=\left\{%
\begin{array}{ll}
    (-1)^{\left\langle\lambda/\nu\right\rangle}\prod\limits_{s\in SC}
    [\Ct(\lambda/\mu)-\Ct(s)]\underset{d\neq\lambda/\mu}{\prod\limits_
{d\in DB}}[\Ct(\lambda/\mu)-\Ct(d)]^{-1},  \\
    \;\;\;\;\;\;\;\;\;\;\;\hbox{if $\lambda/\mu$ is a broken border
strip;} \\
    \\
    0, \;\hbox{otherwise;} \\
\end{array}%
\right.
$$
where $SC$ and $DB$ denote the sets of sharp corners and dull
boxes in $\lambda/\nu$, respectively, and
$\left\langle\lambda/\nu\right\rangle$ is the height of
$\lambda/\nu$.
\end{prop}
\begin{proof} The formula for $\TR(\mu/\nu,\lambda/\nu)$ stated in
the Proposition follows immediately from the relation between
$\TR(\mu/\nu,\lambda/\nu)$ and $X_{\lambda/\nu,\mu/\nu}$, and from
the explicit formula for $X_{\lambda/\nu,\mu/\nu}$ (Theorem
\ref{FORMULAFORXLM}).
\end{proof}

\section{Modules $C'$ and $\Lambda[t]$}\label{SECTION5}
\subsection{The ring $C$}\label{subsection3.1}
Let $S(n)$ be the group of permutations of the set $\{1,2,\ldots
,n\}$. Denote by $C(n)$ the space of (complex valued) functions on
the group $S(n)$ invariant with respect to conjugations by the
elements of $S(n)$. (Thus $C(n)$ is the space of central functions
on $S(n)$).

Let
$$
C=\underset{n\geq 0}\bigoplus\; C(n),
$$
with the understanding that $C(0)=\mathbb{C}$. We embed
$S(m)\times S(n)$ in $S(m+n)$ in such a way that $S(m)$ is the
group of permutations of the first $m$ elements in the set
$\{1,2,\ldots ,m+n\}$, and $S(n)$ is the group of permutations of
the subsequent $n$ elements.

Let $f\in C(m)$, $g\in C(n)$, and we define a bilinear
multiplication $C(m)\times C(n)\rightarrow C(m+n)$ by the formula
\begin{equation}\label{511}
f\cdot g=\Ind_{S(m)\times S(n)}^{S(m+n)}\left(f\times g\right).
\end{equation}
More explicitly, for any $\pi\in S(m+n)$
$$
(f\cdot g)(\pi)=\frac{1}{|S(m)||S(n)|}\sum\limits_{\sigma\in
S(m+n)}(f\times g)(\sigma^{-1}\pi\sigma),
$$
with the understanding that $(f,g)$ is concentrated on $S(m)\times
S(n)$, i.e.
$$
(f\times g)(\sigma^{-1}\pi\sigma)=\left\{%
\begin{array}{ll}
    f(\pi_1)g(\pi_2), & \sigma^{-1}\pi\sigma=\pi_1\pi_2, \pi_1\in S
(m),\pi_2\in S(n); \\
    0, & \hbox{otherwise.} \\
\end{array}%
\right.
$$
Thus $f\cdot g$ is  an element of $C(m+n)$, and with this
multiplication $C$ is a commutative, associative, graded ring with
an identity element.
\subsection{The module $C'$ over $C$.}
Suppose $S(n-1)$ is the subgroup of $S(n)$ realized as the group
of permutations of the set $\{1,2,\ldots ,n-1\}$. Denote by
$C'(n)$ the space of (complex-valued) functions on $S(n)$
invariant with respect to conjugations by elements of the subgroup
$S(n-1)$ of the group $S(n)$. We agree that $S(0)=\{e\}$, and set
$$
C'=\underset{n\geq 1}\bigoplus\; C'(n).
$$
Let us embed $S(m)\times S(n)$ in $S(m+n)$ as in $\S$
\ref{subsection3.1}. Let $f\in C(m)$, and $g\in C'(n)$. We define
a bilinear map
$$
C(m)\times C'(n)\rightarrow C'(m+n)
$$
by the formula
\begin{equation}\label{Eq321}
(f*g)(\pi)=\frac{1}{|S(m)||S(n-1)|}\sum\limits_{\sigma\in
S(m+n-1)}(f\times g)(\sigma^{-1}\pi\sigma),
\end{equation}
where $\pi$ is an arbitrary element of the group $S(m+n)$, and
$(f\times g)$ is concentrated on the subgroup $S(m)\times S(n)$ of
$S(m+n)$. With this bilinear map $C'$ is a module over the ring
$C$.
\subsection{The module $\Lambda[t]$ over $\Lambda$. Isomorphism between
$C'$ and $\Lambda[t]$} As it is well known the ring $C$ is
isomorphic to the ring of symmetric functions $\Lambda$. This
isomorphism is given by the characteristic map, see Macdonald
\cite{macdonald}, I,$\S 7$.

Let $\Lambda[t]$ be the set of polynomials in $t$ whose
coefficients are elements of $\Lambda$. Clearly, $\Lambda[t]$ can
be understood as a module over the ring $\Lambda$. We aim to
construct an isomorphism between $C'$ and $\Lambda[t]$.

Let $w$ be an element of $S(n)$, and suppose $w$ belongs to a
conjugacy class (with respect to the conjugations by $S(n-1)$)
parameterized by the pair $(j,\rho)$. Here $j$ takes values from
$1$ to $n$, and $\rho$ is a partition of $n-j$.

We then define a mapping $ \psi:\; S(n)\rightarrow \Lambda[t]$ as
follows
\begin{equation}\label{330}
\psi(w)=t^{j-1}p_{\rho}
\end{equation}
where $p_{\rho}=p_{\rho_1}p_{\rho_2}\ldots $, and $p_k$ is
$k^{\mbox{th}}$ power sum. Next we define a mapping
$$
\Ch':\; C'\rightarrow \Lambda[t]
$$
in the following way: if $g\in C'(n)$, then
\begin{equation}\label{Eq331}
\Ch'(g)=\frac{1}{(n-1)!}\sum\limits_{w\in S(n)}g(w)\psi(w).
\end{equation}
We note that if $g_{(j,\rho)}$ is the value of $g$ at the
conjugacy class parameterized by $(j,\rho)$, we obtain
\begin{equation}\label{Eq332}
\Ch'(g)=\sum\limits_{j=1}^nt^{j-1}\sum\limits_{\lambda\vdash
n-j}z_{\lambda}^{-1}g_{(j,\rho)}p_{\rho}
\end{equation}
where $z_{\rho}=\prod_{i\geq 1}i^{m_i}m_i!$, $m_i=m_i(\rho)$ is
the number of parts of $\rho$ equal to $i$.

We call $\Ch'(f)$ the characteristic of $f$, and call $\Ch'$ the
characteristic map.
\begin{prop} The characteristic map $\Ch'$ is an isomorphism
between the $C$-module $C'$ and the $\Lambda$-module $\Lambda[t]$.
\end{prop}
\begin{proof}
It is enough to show that for every $m\geq 0$, $n\geq 1$, and for
every $f\in C(m)$, $g\in C'(n)$ the following formula holds:
$$
\Ch'(f*g)=\Ch(f)\cdot\Ch'(g)
$$
where $\Ch(f)$ is defined as in Macdonald \cite{macdonald}, I, $\S
7$.

The left-hand part of this equation is given by formulas
(\ref{Eq321}), (\ref{Eq331}). Using the fact that $\psi(w)$ is
constant on conjugacy classes defined with respect to conjugations
by $S(n-1)$, we rewrite $\Ch'(f*g)$ as follows
$$
\Ch'(f*g)=\frac{1}{|S(m)||S(n-1)|}\sum\limits_{w\in
S(m+n)}(f,g)(w)\psi(w).
$$
Since $(f\times g)$ is concentrated on the subgroup $S(m)\times
S(n)$ of the group $S(m+n)$ we can rewrite the expression above as
a double sum:
$$
\Ch'(f*g)=\frac{1}{|S(m)||S(n-1)|}\sum\limits_{\sigma\in
S(m)}\sum\limits_{\pi\in S(n)}
f(\sigma)g(\pi)\psi(\sigma\cdot\pi).
$$
Suppose $\pi\in S(n)$, $\sigma\in S(m)$, $\pi$ has the cyclic
structure $(j,\mu)$, where $1\geq j\geq n$,  $\mu$ is a partition
of $n-j$, and $\sigma$ has the cyclic structure $\lambda$. Then
the permutation $\sigma\cdot\pi$ has the cyclic structure
$(j,\lambda\cup\mu)$, and we obtain:
$$
\psi(\sigma\cdot\pi)=y^{j-1}\cdot p_{\lambda\cup\mu}=y^{j-1}
p_{\lambda}\cdot p_{\mu}.
$$
Let $f_{\lambda}$ be the value of $f_{\lambda}$ on the
permutations of $S(m)$ with the cyclic structure $\lambda$, and
$g_{(j,\mu)}$ be the value of $g$ on the permutations with the
cyclic structure $(j,\mu)$. Then the expression for $\Ch'(f*g)$
can be rewritten as
$$
\Ch'(f*g)=\left(\sum\limits_{\lambda\vdash
m}z^{-1}_{\lambda}f_{\lambda}p_{\lambda}\right)
\left(\sum\limits_{j=1}^n\sum\limits_{\mu\vdash
n-j}z^{-1}_{\mu}y^{j-1}p_{\mu}g_{j,\mu}\right),
$$
i.e. $\Ch'(f*g)=\Ch(f)\cdot\Ch'(g)$.
\end{proof}
\subsection{Scalar products on $C'$ and $\Lambda[t]$}
Let $f,g\in C'$, say $f=\sum f_n$ and $g=\sum g_n$ with $f_n,
g_n\in C'(n)$. We set
$$
\left\langle f,g\right\rangle=\sum\limits_{n\geq 1}\left\langle
f_n,g_n\right\rangle'_{S(n)}
$$
where
$$
\left\langle
f_n,g_n\right\rangle'_{S(n)}=\frac{1}{(n-1)!}\sum\limits_{w\in
S(n)}f_n(w)g_n(w^{-1}).
$$
This defines a scalar product on $C'$. We now define a scalar
product on $\Lambda[t]$, i.e. a $\mathbb{C}$-valued bilinear form
$\left\langle u,v\right\rangle$, by requiring that the elements of
the basis $\{t^kp_{\lambda}\}$ of $\Lambda[t]$ satisfy the
following orthogonality relation:
\begin{equation}\label{241}
\left\langle
t^kp_{\lambda},t^jp_{\mu}\right\rangle=z_{\lambda}\delta^{k,j}\delta_
{\lambda\mu}.
\end{equation}
\begin{prop}
The characteristic map $\Ch'$ is an isometry, i.e. for any $f$ and
$g$ in $C'(n)$,
\begin{equation}\label{243}
\left\langle\Ch'(f),\Ch'(g)\right\rangle=\left\langle
f,g\right\rangle'_{S(n)}.
\end{equation}
\end{prop}
\begin{proof}
We rewrite the left-hand of equation (\ref{243}) explicitly using
formula (\ref{Eq332}). Applying the orthogonality relation for the
elements of the basis $\{t^kp_{\lambda}\}$ (equation \ref{241}) we
obtain
$$
\left\langle\Ch'(f),\Ch'(g)\right\rangle=\sum\limits_{k=1}^n\sum\limits_
{\lambda\vdash n-k}z_{\lambda}^{-1}f_{(k,\lambda)}g_{(k,\lambda)}
$$
Clearly, the right-hand side of this expression can be rewritten
as $\left\langle f,g\right\rangle'_{S(n)}$.
\end{proof}
\section{The generalized Schur functions}\label{SECTION6}
\subsection{Definition of the generalized Schur functions}
Let $\lambda$ be a partition of $n$, $\mu$ be a partition of
$n-1$. Suppose $\mu$ is obtained from $\lambda$ by removing one
box. Thus we have a pair $(\lambda,\mu\nearrow\lambda)$ of
partitions. Recall that the generalized characters are
parameterized by such pairs.
\begin{defn}
Let $\Gamma^{\lambda,\mu\nearrow\lambda}$ be the generalized
character associated with the Gelfand pair $(S(n)\times
S(n-1),\Diag S(n-1))$. The image $S^{\lambda,\mu\nearrow\lambda}$
of $\Gamma^{\lambda,\mu\nearrow\lambda}$ under the characteristic
map $\Ch'$ is called the generalized Schur function associated
with the Gelfand pair $(S(n)\times S(n-1),\Diag S(n-1))$. Thus
\begin{equation}\label{31}
S^{\lambda,\mu\nearrow\lambda}=\Ch'(\Gamma^
{\lambda,\mu\nearrow\lambda}).
\end{equation}
\end{defn}
\subsection{Orthogonality}
Since the map $\Ch'$ is an isometry, and the generalized
characters satisfy the orthogonality condition
$$
\left\langle\Gamma^{\lambda,\mu\nearrow\lambda},
\Gamma^{\rho,\nu\nearrow\rho}\right\rangle'_{S(n)}=n\;\frac{\dim\mu}
{\dim\lambda}\delta^{\lambda,\rho}\delta^{\nu,\mu},
$$
(where $n=|\lambda|$) it follows that the generalized Schur
functions are orthogonal with respect to the scalar product on
$\Lambda[t]$. Namely,
\begin{equation}\label{32}
\left\langle S^{\lambda,\mu\nearrow\lambda},
S^{\rho,\nu\nearrow\rho}\right\rangle=n\;\frac{\dim\mu}{\dim\lambda}
\delta^{\lambda,\rho}\delta^{\nu,\mu}.
\end{equation}
\begin{prop}
We have
\begin{equation}
\sum\limits_{\lambda,\mu\nearrow\lambda}n^{-1}\frac{\dim\lambda}
{\dim\mu}S^{\lambda,\mu\nearrow\lambda}(t;x_1,x_2,\ldots)
S^{\lambda,\mu\nearrow\lambda}(t;y_1,y_2,\ldots)=\frac{1}{1-t^2}\;\frac
{1}{\prod_{i,j}(1-x_iy_j)}
\end{equation}
where $n=|\lambda|$.
\end{prop}
\begin{proof}
By (\ref{32}), $(S^{\lambda,\mu\nearrow\lambda})$ and
$(n^{-1}\frac{\dim\lambda}{\dim\mu}\cdot
S^{\lambda,\mu\nearrow\lambda})$ are dual bases of $\Lambda[t]$
for the scalar product $\left\langle f,g\right\rangle$ defined by
formula (\ref{241}). $(t^kp_{\lambda})$ and
$(z_{\lambda}^{-1}t^kp_{\lambda})$ are dual bases as well. It
follows (cf. Macdonald \cite{macdonald}, I, $\S 4$, (4.6)) that
$$
\sum\limits_{\lambda,\mu\nearrow\lambda}n^{-1}\frac{\dim\lambda}
{\dim\mu}S^{\lambda,\mu\nearrow\lambda}(t;x_1,x_2,\ldots)
S^{\lambda,\mu\nearrow\lambda}(t;y_1,y_2,\ldots)=
\sum\limits_{k,\lambda}t^{2k}z_{\lambda}^{-1}p_{\lambda}(x_1,x_2,\ldots)
p_{\lambda}(y_1,y_2,\ldots)
$$
$$
=\frac{1}{1-t^2}\;\frac{1}{\prod_{i,j}(1-x_iy_j)}.
$$
\end{proof}
\subsection{Frobenius type formula}
Now we obtain an analogue of the Frobenius formula for the
generalized characters. We have from (\ref{31}) and (\ref{Eq332})
\begin{equation}\label{630}
S^{\lambda,\mu\nearrow\lambda}=\Ch'(\Gamma^{\lambda,\mu\nearrow\lambda})
=\sum\limits_{j=1}^n t^{j-1}\sum\limits_{\rho\vdash
n-j}z_{\rho}^{-1}\Gamma^{\lambda,\mu\nearrow\lambda}_{(j,\rho)}\cdot
p_{\rho}
\end{equation}
where $\Gamma^{\lambda,\mu\nearrow\lambda}_{(j,\rho)}$ is the
value of $\Gamma^{\lambda,\mu\nearrow\lambda}$ at elements of
cycle-type $(j,\rho)$. Hence
$$
\Gamma^{\lambda,\mu\nearrow\lambda}_{(j,\rho)}=\left\langle
S^{\lambda,\mu\nearrow\lambda}, t^{j-1}p_{\rho}\right\rangle,
$$
i.e. the transition matrix between the bases
$(S^{\lambda,\mu\nearrow\lambda})$ and $(t^{j-1}p_{\rho})$ is the
table of the generalized characters. Thus
 \begin{equation}\label{34}
t^{j-1}p_{\rho}=\sum\limits_{\lambda,\mu\nearrow\lambda}\Gamma^
{\lambda,\mu\nearrow\lambda}_{(j,\rho)}
S^{\lambda,\mu\nearrow\lambda}
\end{equation}
where the sum is over all partitions $\lambda$ of $n$, and
partitions $\mu$ of $n-1$ such that $\mu\nearrow\lambda$.
\subsection{Formula for the generalized Schur functions}
The Murnaghan-Nakayama type rule for the generalized characters
enables us to give an explicit formula for the generalized Schur
functions.
\begin{thm}
Let $\lambda$ be a partition of $n$, $\mu$ be a partition of
$n-1$, and $\mu$ is obtained from $\lambda$ by removing one box.
Then the following formula holds
\begin{equation}\label{641}
S^{\lambda,\mu\nearrow\lambda}(t;x_1,x_2,\ldots)=\sum\limits_
{\nu\subseteq\mu}\varphi_{\mu/\nu,\lambda/\nu}t^{|\lambda/\nu|-1}
s_{\nu}(x_1,x_2,\ldots)
\end{equation}
where $s_{\nu}$ is the Schur symmetric function associated with
the partition $\nu$, and $\varphi_{\mu/\nu,\lambda/\nu}$ is a
combinatorial coefficient associated with the skew Young diagram
$\lambda/\nu$. This combinatorial coefficient is given explicitly
in Theorem \ref{MAINTHEOREM}.
\end{thm}
\begin{proof}
Let $\theta$ be a permutation from $S(n)$ which is in the standard
form
$$
\theta=(1,2,\ldots ,i_1)(i_1+1,\ldots ,i_2)\ldots
(i_{m-1}+1,\ldots ,n)
$$
Suppose $\theta$ has a cyclic structure $(k,\rho)$ (with respect
to the conjugations by $S(n-1)$). Then the last cycle of $\theta$
(those one which includes $n$) has length $k$. Denote this cycle
by $\sigma$, and denote by $\pi$ the permutation of the remaining
$n-k$ numbers. Thus $\theta=\pi\cdot\sigma$, and we have
$$
\Gamma^{\lambda,\mu\nearrow\lambda}(\pi\cdot\sigma)=\Gamma^
{\lambda,\mu\nearrow\lambda}_{(k,\rho)}.
$$
Theorem \ref{MAINTHEOREM} gives the following formula for the
generalized character
$\Gamma^{\lambda,\mu\nearrow\lambda}_{(k,\rho)}$:
$$
\Gamma^{\lambda,\mu\nearrow\lambda}_{(k,\rho)}=\underset{\nu\vdash
n-k}{\sum\limits_{\nu\subseteq\mu
}}\varphi_{\mu/\nu,\lambda/\nu}\chi^{\nu}_{\rho}.
$$
Let us insert this expression into formula (\ref{630}). We find
\begin{equation}
\begin{split}
S^{\lambda,\mu\nearrow\lambda}&=\sum\limits_{k=1}^n
t^{k-1}\sum\limits_{\rho\vdash
n-k}z_{\rho}^{-1}\Gamma^{\lambda,\mu\nearrow\lambda}_{(k,\rho)}\cdot
p_{\rho}\\
&=\sum\limits_{k=1}^n t^{k-1}\sum\limits_{\rho\vdash
n-k}z_{\rho}^{-1}\left(\underset{\nu\vdash
n-k}{\sum\limits_{\nu\subseteq\mu
}}\varphi_{\mu/\nu,\lambda/\nu}\chi^{\nu}_{\rho}\right)
p_{\rho}\\
&=\sum\limits_{k=1}^n t^{k-1}\underset{\nu\vdash
n-k}{\sum\limits_{\nu\subseteq\mu
}}\varphi_{\mu/\nu,\lambda/\nu}\sum\limits_{\rho\vdash
n-k}z_{\rho}^{-1}\chi^{\nu}_{\rho}
p_{\rho}\\
&=\sum\limits_{k=1}^n t^{k-1}\underset{\nu\vdash
n-k}{\sum\limits_{\nu\subseteq\mu
}}\varphi_{\mu/\nu,\lambda/\nu}s_{\nu},\\
\nonumber
\end{split}
\end{equation}
which is obviously equivalent to the formula in the statement of
the Theorem. (We have used the well-known formula
$$
s_{\nu}=\sum\limits_{\rho\vdash n-k}z_{\rho}^{-1}\chi^{\nu}_{\rho}
p_{\rho}
$$
in the last equation.)
\end{proof}
\subsection{The Jacobi-Trudi type formula}
The Schur symmetric functions can be expressed as polynomials in
the complete symmetric functions. The formula is
\begin{equation}\label{JacobiTrudiFormula}
s_{\lambda}=\det\left(h_{\lambda_i-i+j}\right)_{1\leq i,j\leq n},
\end{equation}
where $n\geq l(\lambda)$, $l(\lambda)$ denotes the number of rows
of partition $\lambda$. We are looking for an analogue of this
formula for the generalized Schur functions.

Let $1_{n}$ be the identity character of $S(n)$. Then
\begin{equation}\label{651}
\Ch(1_n)=\sum\limits_{|\rho|=n}z_{\rho}^{-1}p_{\rho}=h_n.
\end{equation}
If $\lambda=(\lambda_1,\lambda_2,\ldots )$ is any partition of
$n$, let $1_{\lambda}$ denote $1_{\lambda_1}\cdot
1_{\lambda_2}\ldots$ with the multiplication defined by formula
(\ref{511}). Then $1_{\lambda}$ is the character of $S(n)$ induced
by the identity character of $S_{\lambda}=S_{\lambda_1}\times
S_{\lambda_2}\ldots$, and we have $\Ch(1_{\lambda})=h_{\lambda}$.
Moreover, formula (\ref{651}) enables to rewrite the Jacobi-Trudi
formula, equation (\ref{JacobiTrudiFormula}), as an  expression
for the irreducible characters of $S(n)$ in terms of the induced
characters:
$$
\chi^{\lambda}=\det\left(1_{\lambda_i-i+j}\right)_{1\leq i,j\leq
n}.
$$
Denote by $1_k'$ the identity character of $S(k)$  considered  as
an element of $C'(k)$.
\begin{prop}
For every $k\geq 1$, and every partition $\mu=(\mu_1,\ldots
,\mu_m)$ we obtain
\begin{equation}
\Ch'(1_{\mu_1}*\ldots
* 1_{\mu_m}*1_k')=h_{\mu_1}\cdot\ldots\cdot h_{\mu_m}\cdot h_{k}'
\end{equation}
where
\begin{equation}\label{654}
h_k'=\sum\limits_{j=1}^kt^{j-1}h_{k-j}.
\end{equation}
\end{prop}
\begin{proof}
We have
\begin{equation}
\Ch'(1_k')=\sum\limits_{j=1}^kt^{j-1}\sum\limits_{\lambda\vdash
k-j}z_{\lambda}^{-1}p_{\lambda}=\sum\limits_{j=1}^kt^{j-1}h_{k-j}=h_k'
\end{equation}
where  we have used formula (\ref{Eq332}) in the first equation.
Since $\Ch'(f*g)=\Ch(f)\Ch'(g)$ for $f\in C(m)$ and $g\in C'(k)$,
the formula in the statement of the Proposition holds.
\end{proof}
Now note that equation (\ref{654}) implies
$$
t^{n-1}=(-1)^{n-1}\left|\begin{array}{ccccc}
  h_1 & h_2 & h_3 & \ldots & h_n' \\
  1 & h_1 & h_2 & \ldots & h_{n-1}' \\
  0 & 1 & h_1 & \ldots & h_{n-2}' \\
  \vdots & \ddots &  &  &  \\
  0 & 0 & \ldots & 1 & h_1' \\
\end{array}\right|.
$$
If we insert this expression into formula (\ref{641}) and rewrite
$s_{\nu}$ as the Jacobi-Trudi determinant we obtain a Jacobi-Trudi
like representation of the generalized Schur functions:
\begin{equation}\label{655}
\begin{split}
S^{\lambda,\mu\nearrow\lambda}&=\sum\limits_{\nu\subseteq\mu}\varphi_
{\lambda/\nu}(-1)^{|\lambda/\nu|-1}
\Biggl\{\left|\begin{array}{ccccc}
  h_1 & h_2 & \ldots & h_{|\lambda/\nu|-1} & h_{|\lambda/\nu|}' \\
  1 & h_1 & h_2 & h_{|\lambda/\nu|-2} & h_{|\lambda/\nu|-1}' \\
  0 & 1 & \ldots & h_{|\lambda/\nu|-3} & h_{|\lambda/\nu|-2}' \\
  \vdots & \ddots &  &  &  \\
  0 & 0 & \ldots & 1 & h_1' \\
\end{array}\right|\\
&\times\left|\begin{array}{cccc}
  h_{\nu_1} & h_{\nu_1+1} & \ldots & h_{\nu_1+l(\nu)-1} \\
  h_{\nu_2-1} & h_{\nu_2} & \ldots & h_{\nu_2+l(\nu)-2} \\
  \vdots & \vdots &  &  \\
  h_{\nu_l-l(\nu)+1} & h_{\nu_l-l(\nu)+2} & \ldots  & h_{\nu_{l}} \\
\end{array}\right|\Biggr\}.
\end{split}
\end{equation}
This is an analogue of the Jacobi-Trudi formula for the
generalized Schur functions.\\


\begin{thebibliography}{0000}
\bibitem[AkMiz]{akazawa} H. Akazawa and H. Mizukawa.
 \textit{Orthogonal polynomials arising from
the wreath products of a dihedral group with a symmetric group.}
J. Combin. Theory Ser. A \textbf{104} (2003) 371--380.


\bibitem[Ban]{bannai} E. Bannai, T. Ito. Algebraic Combinatorics I.
Association Schemes, The Benjamin/Cummings Publishing Co., CA,
1984.
\bibitem[Bren]{bender} M. Brender.
\textit{Spherical Functions on the Symmetric Group.} J. of Algebra
\textbf{42} (1976), 302--314.

\bibitem[Gal]{gallagher}
P. Gallagher. \textit{Functional equation for spherical functions
on finite groups}, Math. Z. \textbf{141} (1975), 77--81.

\bibitem[GW]{goodman} R. Goodman and N. Wallach. Representations
and invariants of the classical groups. (Encycl. Math. Appl.
\textbf{68}) Cambridge: Cambridge U. Press 1998.

\bibitem[Gr]{greene} C. Greene. \textit{A rational-function identity
related to the Murnaghan-Nakayama formula for the characters of
$S_n$.} J. Algebraic Combinatorics \textbf{1} (1992), 235--255.

\bibitem[Ju]{jucys} A. Jucys. \textit{Symmetric polynomials and the center
of the symmetric group ring}. Reports Math. Phys. \textbf{5}
(1974), 107-112.

\bibitem[HalRam]{halverson} T. Halverson and A. Ram.
\textit{Murnaghan-Nakayama rules for characters of Iwahori-Hecke
algebras of classical type.} Trans. Amer. Math. Soc. \textbf{348}
(1996), 3967-3995.

\bibitem[Knop]{knop}
F. Knop. \textit{Semisymmetric polynomials and the invariant
theory of matrix vector pairs.}  Represent. Theory \textbf{5}
(2001), 224--266.

\bibitem[Mac]{macdonald}
 I., G. Macdonald. Symmetric functions and Hall polynomials, 2nd
edition, Oxford University Press, 1995.

\bibitem[Miz]{mizukawa}
H. Mizukawa. \textit{Zonal spherical functions on the complex
reflection groups and (n+1,m+1)-hypergeometric functions.} Adv. in
Math. \textbf{184} (2004) 1-17.



\bibitem[MizTan]{mizukawa1} H. Mizukawa and H. Tanaka.  \textit{(n+1, m+1)- hypergeometric functions
associated to character algebras.} Proc. Amer.Math. Soc.
\textbf{132} (2004) 2613-2618.

\bibitem[Mur]{murphy} G. Murphy. \textit{A new construction of Young's
seminormal representation of the symmetric group.} J. Algebra
\textbf{69} (1981), 287-291.

\bibitem[Ok1]{okounkov1} \textit{Thoma's theorem and representations of
infinite bisymmetric group}. Func. Anal. Appl. \textbf{28} (1994),
no. 2,101-107.

\bibitem[Ok2]{okounkov2} \textit{On representations of infinite symmetric
group.}  Zap. Nauchn. Sem. S.-Peterburg. Otdel. Mat. Inst.
Steklov. (POMI) \textbf{240} (1997), Teor. Predst. Din. Sist.
Komb. i Algoritm. Metody. \textbf{2}, 166--228, 294; translation
in J. Math. Sci. (New York) \textbf{96} (1999), no. 5, 3550--3589.


\bibitem[OkVer]{okounkov}
A. Okounkov and A. Vershik. \textit{A new approach to
representation theory of symmetric groups.} Selecta Mathematica
\textbf{1} (1996), 581-605.

\bibitem[Olsh]{olshanski}
G. I. Olshanski. \textit{Unitary representations of $(G,K)$-pairs
connected with the infinite symmetric group $S(\infty)$.}
Leningrad Math. J. \textbf{1} (1990), 985--1014.

\bibitem[Olsh1]{olshanski1}
G. I. Olshanski. \textit{Extension of the algebra $U(g)$ for
infinite-dimensional classical Lie algebras $g$ and the Yangians
$Y(gl(m))$}. Soviet Math. Dokl. \textbf{36} (1988), 569-573.

\bibitem[Roich]{roichman}
Y. Roichman. Characters of the symmetric groups: formulas,
estimates and applications. Emerging applications of number theory
(Minneapolis, MN, 1996), 525--545, IMA Vol. Math. Appl., 109,
Springer, New York, 1999.


\bibitem[Ruth]{rutherford}
D., E. Rutherford. Substitutional Analysis,. University press,
Edinburg, 1948.

\bibitem[S]{sagan} B. Sagan. The symmetric group. Representations,
combinatorial algorithms, and symmetric functions. Second edition.
Graduate Texts in mathematics, 203. Springer-Verlag, New-York,
2001.
\bibitem[St]{stanley} R. P. Stanley. Enumerative combinatorics,
Vol. 2, Cambridge University Press, San Diego, 1991.

\bibitem[Trav]{travis} D. Travis. \textit{Spherical Functions of Finite
Groups.} J. of Algebra \textbf{29} (1974), 65--76.

\bibitem[Young]{young} A. Young. The collected papers of Alfred
Young. University of Toronto Press, Toronto, 1977.
\end{thebibliography}
\end{document}